\documentclass[12pt]{amsart}
\usepackage{amsmath,amssymb,amscd,amsthm,latexsym}
\usepackage{euscript}
\DeclareMathAlphabet{\EuRm}{U}{eur}{m}{n}
\SetMathAlphabet{\EuRm}{bold}{U}{eur}{b}{n}
\begin{document} 
\baselineskip 17pt plus2pt minus1pt
\parskip 3pt plus 1pt minus .5pt

%
%
\swapnumbers
\newtheorem{thm}{Theorem}[section]
\newtheorem{lemma}[thm]{Lemma}
\newtheorem{prop}[thm]{Proposition}
\newtheorem{cor}[thm]{Corollary}
\theoremstyle{definition}
\newtheorem{defn}[thm]{Definition}
\newtheorem{example}[thm]{Example}
\newtheorem{examples}[thm]{Examples}
\newtheorem{claim}[thm]{Claim}
\newtheorem{con}[thm]{Convention}
\newtheorem{obs}[thm]{Observation}

\theoremstyle{remark}
\newtheorem{assume}[thm]{Assumption}
\newtheorem{remark}[thm]{Remark}
\newtheorem{note}[thm]{Note}
\newtheorem{notation}[thm]{Notation}
\newtheorem{aside}[thm]{Aside}
\newtheorem{ack}[thm]{Acknowledgements}
\numberwithin{equation}{section}
\numberwithin{figure}{section}
%
%
\def\sect{\setcounter{thm}{0}\section}
%
%
\newcommand{\xra}[1]{\xrightarrow{#1}}
\newcommand{\xla}[1]{\xleftarrow{#1}}
\newcommand{\xsim}{\xrightarrow{\sim}}
\newcommand{\hra}{\hookrightarrow}
\newcommand{\epic}{\to\hspace{-5 mm}\to}
\newcommand{\adj}[2]{\substack{{#1}\\ \rightleftharpoons \\ {#2}}}
\newcommand{\csub}[1]{\hspace{1mm}\circ_{#1}\hspace{1mm}}
\def\ccsub#1{\circ_{#1}}
\newcommand{\tripless}[3]{\def\arraystretch{.5}
\begin{array}{c}\mbox{\scriptsize $\scriptstyle{#1}$}\\
\mbox{\scriptsize $\scriptstyle{#2}$}\\
\mbox{\scriptsize $\scriptstyle{#3}$}\end{array}\def\arraystretch{1}}
\newcommand{\DEF}{:=}
\newcommand{\EQUIV}{\Leftrightarrow}
\newcommand{\hsp}{\hspace{10 mm}}
\newcommand{\hs}{\hspace{5 mm}}
\newcommand{\hsm}{\hspace{3 mm}}
\newcommand{\vsm}{\vspace{2 mm}}
\newcommand{\rest}[1]{\lvert_{#1}}
\newcommand{\lra}[1]{\langle{#1}\rangle}
\newcommand{\llrr}[1]{\langle\!\langle{#1}\rangle\!\rangle}
\newcommand{\llrra}[1]{\overline{\langle\!\langle{#1}\rangle\!\rangle}}
\newcommand{\act}[1]{\circ{#1}}
%
%

\newcommand{\squar}[8]{\setlength{\unitlength}{.7cm}
\begin{picture}(4,3.7)(0,0)
\thinlines
\put(0,3){\makebox(0,0){$#1$}}
\put(5,3){\makebox(0,0){$#2$}}
\put(0,0){\makebox(0,0){$#3$}}
\put(5,0){\makebox(0,0){$#4$}}
\put(-.5,1.5){\makebox(0,0)[r]{$#6$}}
\put(5.5,1.5){\makebox(0,0)[l]{$#7$}}
\put(2.5,0.5){\makebox(0,0)[b]{$#8$}}
\put(2.5,3.5){\makebox(0,0)[b]{$#5$}}
\put(1,0){\vector(1,0){3}}
\put(1,3){\vector(1,0){3}}
\put(0,2.5){\vector(0,-1){2}}
\put(5,2.5){\vector(0,-1){2}}
\end{picture}}
%
%
\def\vlra{{\hbox{$-\hskip-1mm-\hskip-2mm\longrightarrow$}}}
\newcommand{\A}{{\EuScript A}}
\newcommand{\BA}{C\A}
\newcommand{\bA}{\bar{\A}}
\newcommand{\B}{{\EuScript B}}
\newcommand{\Ad}{\A_{\bullet}}
\newcommand{\Ass}{{\EuScript Ass}}
\newcommand{\C}{{\mathcal C}}
\newcommand{\Cu}{{\EuScript Cu}}
\newcommand{\D}{{\mathcal D}}
\newcommand{\F}{{\EuScript F}}
\newcommand{\II}{{\mathcal I}}
\newcommand{\OO}{{\EuScript O}}
\newcommand{\PP}{{\mathcal P}}
\newcommand{\bPP}{b\PP}
\newcommand{\PG}{\PP_{\Gamma}}
\newcommand{\Perm}{{\EuScript Perm}}
\newcommand{\Ss}{{\mathcal S}}
\newcommand{\TT}{{\mathcal T}}
\newcommand{\Ta}{\TT_{\ast}}
\newcommand{\hoT}{\ho\Ta}
\newcommand{\U}{{\mathcal U}}
\newcommand{\V}{{\mathcal V}}
\newcommand{\WW}[1]{W{#1}}
\newcommand{\WG}{\WW{\Gamma}}
\newcommand{\WsG}{W_{s}\Gamma}
%
%
\newcommand{\N}{{\mathbb N}}
\newcommand{\Q}{{\mathbb Q}}
\newcommand{\R}{{\mathbb R}}
\newcommand{\Z}{{\mathbb Z}}
%
%
\newcommand{\pis}{\pi_{\ast}}
%
%
\newcommand{\sk}[1]{\operatorname{sk}_{#1}}
\newcommand{\Dom}{\operatorname{Dom}}
\newcommand{\End}{\operatorname{End}}
\newcommand{\Ext}{\operatorname{Ext}}
\newcommand{\fib}{\operatorname{fib}}
\newcommand{\fin}{\operatorname{fin}}
\newcommand{\ho}{\operatorname{ho}}
\newcommand{\holim}{\operatorname{holim}}
\newcommand{\Hom}{\operatorname{Hom}}
\newcommand{\Id}{\operatorname{Id}}
\newcommand{\init}{\operatorname{init}}
\newcommand{\vf}{v_{\fin}}
\newcommand{\vfi}[1]{v{#1}_{\fin}}
\newcommand{\vi}{v_{\init}}
\newcommand{\vin}[1]{v{#1}_{\init}}
\newcommand{\map}{\operatorname{Map}}
\newcommand{\Mor}{\operatorname{Mor}}
\newcommand{\Obj}{\operatorname{Obj}\,}
\newcommand{\op}{\operatorname{op}}
\newcommand{\pt}{\operatorname{pt}}
\newcommand{\PROP}{{\sc{PROP}}}
\newcommand{\vdim}{\operatorname{vdim}}
%
%
\newcommand{\BW}[1]{b\WW{#1}}
\newcommand{\BG}{\BW{\Gamma}}
\newcommand{\hBG}{\widehat{\BG}}
\newcommand{\BsG}{bW_{s}\Gamma}
\newcommand{\bD}[1]{\mathbf{D}^{#1}}
\newcommand{\Del}{\text{\boldmath $\Delta$}}
\newcommand{\De}[1]{\Delta[{#1}]}
\newcommand{\bh}{\mathbf{h}}
\newcommand{\bk}{\text{\boldmath $[k]$}}
\newcommand{\bkm}{\text{\boldmath $[k$-$1]$}}
\newcommand{\bl}{\text{\boldmath $[l]$}}
\newcommand{\bm}{\text{\boldmath $[-1]$}}
\newcommand{\bn}{\text{\boldmath $[n]$}}
\newcommand{\bnm}{\text{\boldmath $[n$-$1]$}}
\newcommand{\bnp}{\text{\boldmath $[n$+$1]$}}
\newcommand{\bo}{\text{\boldmath $[1]$}}
\newcommand{\bt}{\text{\boldmath $[2]$}}
\newcommand{\bz}{\text{\boldmath $[0]$}}
\newcommand{\Pe}[1]{P\!e_{#1}}
\newcommand{\bS}[1]{\mathbf{S}^{#1}}
\newcommand{\bW}{\bar{W}}
\newcommand{\wG}{\PG}
\newcommand{\wcG}{\PP^{c}_{\Gamma}}
\newcommand{\wsG}{\PP^{s}_{\Gamma}}
\newcommand{\X}{\mathbf{X}}
\newcommand{\Xd}{X_{\bullet}}
\newcommand{\XAd}{X[\A]_{\bullet}}
\newcommand{\Y}{\mathbf{Y}}
\newcommand{\fimax}{{\phi_{\max}}}
%
%
\newcommand{\david}{\noindent\textbf{David's Comment:}\\ }
\newcommand{\edavid}{\\ \textbf{End of David's Comment.}\\ }
\newcommand{\dave}{\noindent\textbf{\ [David's change: \ }}
\newcommand{\edave}{\textbf{\ End]}\  }
%
%
\title{Higher homotopy operations}
\author{David Blanc and Martin Markl}
\thanks{Second author partially supported by grant GA AV \v{C}R 1019203}
\address{Dept.\ of Mathematics, Univ. of Haifa, 31905 Haifa, Israel}
\address{Mathematical Inst.\ of the Academy, \v{Z}itn\'{a} 25, 115 67
Prague 1, Czech Republic}
\email{blanc@math.haifa.ac.il\\ markl@math.cas.cz}
\date{June 5, 2001; revised June 12, 2002}
\begin{abstract}
We provide a general definition of higher homotopy operations, encompassing 
most known cases, including higher Massey and Whitehead products, and long 
Toda brackets. These operations are defined in terms of the $W$-construction
of Boardman and Vogt, applied to the appropriate diagram category; we
also  show how some classical families of polyhedra (including
simplices, cubes, associahedra, and permutahedra) arise in this way.
\end{abstract}
\maketitle
%
%

\sect{Introduction}
\label{cint}

Higher homotopy operations have a long history, starting with Toda
brackets, Massey products, and Adem's secondary cohomology operations. 
Secondary homotopy and cohomology operations have been exploited with 
great effect \ -- \ for example, in \cite{AdHI,BJMahT,MPetS,PSteS} \ -- \ 
but beyond attempts by Spanier (see \S \ref{shist} below),
there has been no systematic approach to higher \emph{homotopy}
operations in general, though they have appeared sporadically in the 
literature (e.g., in \cite{GWalkL}), and higher \emph{cohomology} operations
have been studied by Maunder and others (\emph{ibid.}).
However, there has been a certain revival of interest 
in higher operations in algebraic contexts (see, for instance, 
\cite{RetaL,AlldR,TanrH}), which may perhaps justify this attempt to
set up the foundations of the subject on a new basis.

We define higher order homotopy operations as the obstruction to
making a homotopy-commutative diagram of spaces \ 
$\A:\Gamma\to\hoT$ \ \emph{strictly} commutative \ (where $\Gamma$
is a certain finite directed category which we term a ``lattice'').
This obstruction is defined in terms of the well-known $W$-construction
of Boardman and Vogt, which takes a particularly convenient form for
lattices, and yields some interesting families of polyhedra as a side
benefit (see Section \ref{cfp} below).

It should be emphasized that our goal here is
to \emph{define} the concept of a higher homotopy operation, rather
than to describe such an explicit obstruction theory for rectifying
diagrams, as in \cite{DKSmH}. Giving an explicit dictionary for
translating between the cohomology obstructions provided by Dwyer, Kan
and Smith (or Spanier's approach, described below), and
our description in terms of higher operations appears to be a
difficult, though interesting, question.

This different viewpoint allows us to sweep the rather messy description of
coherent vanishing of the lower operations under the carpet, giving
the necessary and sufficient conditions in order that a higher
order operation be defined directly in terms of an appropriate map \ 
$\A(\vi)\rtimes b\PP\to\A (\vf)$ \ (see \S \ref{dhho} below).

It should be observed that the theory we describe here is
still somewhat ad hoc; our goal is to present a uniform treatment of
the main known examples, without trying to obtain the most general definition 
possible. Philosophically, higher homotopy operations are connected
with higher homotopies, which arise when one tries to lift a
commutative diagram in the homotopy category to topological spaces 
(they are thus related in principle to operads, although it seems that
the relation has never been made explicit). Therefore, a more general 
definition would perhaps require a satisfactory homotopy theory of
$n$-categories, which is not yet available.

Note that, in practice, higher operations often appear in an algebraic form,
as differentials in spectral sequences (for example, \cite[Ch.\ 2]{AdHI}, 
or \cite[Prop.\ 4.2.5]{BlaO}), as \ $\Ext$ \ classes 
(cf.\ \cite[Ch.\ 16, 3]{MargS}), and so on; these often serve as an efficient 
means of \emph{computing} such operations (e.g., \cite[\S 6]{BlaH}).
But one should think of the operation itself as the intrinsic 
homotopy-theoretic fact, which may manifest itself in different
(seemingly unrelated) algebraic guises.

\subsection{Notation}
\label{snot}\stepcounter{thm}

The category of compactly generated topological spaces is denoted by $\TT$,  
and that of pointed connected compactly generated spaces by \ $\Ta$. \ Their
homotopy categories are denoted by \ $\ho\TT$ \ and \ $\hoT$ \
respectively. The category of simplical sets will be denoted by \ $\Ss$.

Let \ $\N_{+}$ \ denote the category of finite sets \ 
$\bn\DEF\{0,1,2,\dotsc,n\}$ \ ($n=-1,0,1,\dotsc$, \ 
where \ $\bm\DEF\emptyset$), \ with order-preserving monomorphisms as maps. 
The morphisms are generated by  the inclusions \ 
$d^{i}_{n}:\bn=\{0,1,2,\dotsc,n\}\cong
     \{0,1,\dotsc,\hat{i},\dotsc,n+1\}\hra\bnp$ \ ($i=0,1,\dotsc,n+1$). \ 
We denote by \ $\N$ \ the full subcategory of non-empty 
finite sets in \ $\N_{+}$ \ (i.e., omit \ $\bm$). \ A functor \ 
$T:\N^{\op}\to\C$ \ is called a $\Delta$-\emph{simplicial} 
(or: restricted simplicial) \emph{object} over $\C$; this is
just a simplicial object without degeneracies.

\subsection{Other approaches to higher operations}
\label{shist}\stepcounter{thm}

Toda's definition of what we now call Toda brackets in \cite{TodG} 
(see \cite[Ch.\ I]{TodC} and Example \ref{etb} below) was the first example of 
a secondary homotopy operation \emph{stricto sensu}, although 
Adem's secondary cohomology operations (see \cite{AdemI}), and Massey's 
triple products in cohomology (see \cite{MassN} and Sec.\ \ref{chmp} below) 
appeared at about the same time.  

In all three cases there was no attempt to provide a theoretical
framework for such operations; it was Adams, in \cite[Ch.\ 3]{AdHI}, 
who first produced a general definition of secondary cohomology 
operations, based on ``universal examples'' \ -- \ cohomology 
classes in the fiber of a map between mod $p$ generalized
Eilenberg-Mac Lane objects (GEMs) \ -- \ and
explained how the so-called \emph{stable} secondary operations
correspond to relations in the Steenrod algebra. This approach was
generalized by Maunder in his thesis to $n$-th order cohomology operations, 
using chain complexes (or towers) of length $n$; he also
gave necessary and sufficient conditions in order for a higher-order
operation to be defined, in terms of certain ``pyramids'' of
lower order operations (see \cite{MaunC}). 

Later, Holtzman (in \cite{HoltH}) described an alternative version of $n$-th
order cohomology operations based on the $p$-divisibility of certain
``pseudo primary operations.''  Finally, Klaus (in \cite{KlauT})
defines unstable higher order operations in terms of (not necessarily 
linear) natural transformations on cochains \ -- \ generalizing the approach
of Kristensen in \cite{KristS,KKrisS}. 

However, in all these approaches we do not see the combinatorial
underpinning (i.e., the lattice $\Gamma$) of our more general definition, 
essentially because the homotopy category of cochain complexes (or of
GEMs) is additive. 

In \cite{SpanH}, Spanier gave a general theory of higher order
operations (extending the definition of secondary operations given
in \cite{SpanS}), somewhat similar in spirit to the approach we propose here:
an \ $(n+1)$-st order operation is defined as a set of cohomology classes \ 
$\OO^{n}\subset H^{n}(K;\Gamma_{n})$, \ where $K$ is a simplicial complex 
(corresponding essentially to our \ $b\PP$ \ -- \ see \S \ref{dbas} below), 
and the coefficients \ $\Gamma_{n}$ \ are a \emph{stack}
(``cosheaf of groups'') on $K$, defined \ 
$$
\Gamma_{n}(\sigma) \DEF \pi_{n}\Phi(\sigma)\hsp\text{for}\hs \sigma\in K,
$$
\noindent where \ $\Phi(\sigma)$ \ is the topological space assigned to the
simplex \ $\sigma\in K$ \ by a given \emph{carrier} (``cosheaf of
spaces''), and $\pi_{n}$ is as usual its $n$-th homotopy group.

However, Spanier does not explain how $\Phi$ and $K$ may be constructed
from a diagram \ $\A:\Gamma\to\hoT$, \ so the interpretation of
the higher operations as obstructions to rectification is obscured.

\begin{ack}\label{ack}\stepcounter{subsection}
We would like to thank Rainer Vogt for explaining some of the secrets of 
\cite{BVogHI} to the second author.  \ 
The first author would like to thank the Department of Mathematics at
Northwestern University for their hospitality during the period when
this paper was completed. Thanks are also due to the referee for his
helpful comments and suggestions.
\end{ack}

%
%
\sect{Lattices and the $W$-construction}
\label{clwc}

Motivated by the examples of \cite{BlaHH,BlaHO}, in Section~\ref{chho}
we are going to define higher homotopy operations for pointed topological 
spaces. This section contains some auxiliary material for definitions in
Section~\ref{chho}. 

\begin{defn}\label{dlatt}\stepcounter{subsection}  
For the purposes of this paper, a \emph{lattice} will be a finite directed 
category $\Gamma$ with weak initial object \
$\vi=\vi(\Gamma)$ \ and weak terminal object \ $\vf=\vf(\Gamma)$. \ 
This means simply that for every object (or \emph{node}) $w$ of
$\Gamma$, there is at least one map (or \emph{arrow}) from \ $\vi$ \ to $w$, 
and at least one from $w$ to \ $\vf$ \ -- \ and that moreover there is
a unique \emph{maximal} arrow \ $\fimax:\vi\to\vf$. \ We will always
assume that $\Gamma$ is \emph{non-trivial} \ -- \ that is, there
exists a node $w$ of $\Gamma$ with \ $\vi\neq w\neq\vf$. 

We shall denote the finite set of objects (or nodes) of
$\Gamma$ by $V$. \ A composable sequence of $k$ arrows in $\Gamma$, \  
$\Phi =\lra{v_{k}\xra{f_{k}} v_{k-1}\to\dotsb\to v_{1}\xra{f_{1}} v_{0}}$, \ 
will be called a $k$-\emph{chain}. We sometimes denote it more briefly by \ 
$f_{1}f_{2}\dotsb f_{k}$ \ (note the reversed order).  If it starts at \ 
$\vi$ \ and ends at \ $\vf$, \ we say $\Phi$ is a \emph{maximal} chain.  
\end{defn}

\begin{example}\label{elat}\stepcounter{subsection}  
For each \ $n\geq 1$ \ we denote by \ $L_{n}$ \ the lattice generated
by $n$ composable morphisms 
$$
\vi \xra{f_{1}} v_{1} \xra{f_{2}} v_{2} \dotsb v_{n-1} \xra{f_{n}} \vf
$$
\noindent (and no further relations).
\end{example}

\begin{remark}\label{runit}\stepcounter{subsection}
Because our category $\Gamma$ is directed, and thus has no non-trivial loops,
we shall treat it as a \emph{non-unital} category \ -- \ i.e., we omit the 
identity morphisms, so \ $\Hom_{\Gamma}(v,v)=\emptyset$, \ and we can think of 
$\Gamma$ simply as a directed graph with \emph{commutation relations} 
(which are completely determined by a set of distinguished ``commutative 
triangles'' in $\Gamma$, with sides \ $f$, $g$, \ and \ $f\circ g$).
\end{remark}

\subsection{The $W$-construction}
\label{sbc}\stepcounter{thm}
A lattice $\Gamma$, being a special kind of a category, is in particular a
\emph{colored operad} with a set of colors $V$ (cf.\ \cite[\S 2]{BoaH}), and 
as such it has a ``resolution,'' denoted by \ $\WG$ \ and called the 
\emph{bar construction} on $\Gamma$ by Boardman and Vogt 
(cf.\ \cite[III, \S 1]{BVogHI}). We would like to think of it as a 
cofibrant replacement in a hypothetical model category structure on the 
category of all (small) categories. 

\begin{defn}\label{dwg}\stepcounter{subsection}  
Let \ $\Gamma_{n+1}(u,v)$ \ be the discrete space of all \ $(n+1)$-chains
$$
\lra{u = v_{n+1}\xra{f_{n+1}} v_{n}\to\dotsb\to v_{1}\xra{f_{1}} v_{0}= v}.
$$
\noindent $I : = [0,1]$ \ denotes as usual the unit interval. 

For \ $u,v \in V$, \ let
$$
\WG(u,v) := \bigsqcup_{n \geq 0} \Gamma_{n+1}(u,v) \times I^{n}/ \sim,
$$
with the relation $\sim$ defined as follows: if
$$
f_{1}\csub{t_{1}}f_{2}\csub{t_{2}}\dotsb\csub{t_{n-1}}f_{n}\csub{t_{n}}f_{n+1}
$$
\noindent denotes the point 
$$
\lra{u\xra{f_{n+1}} v_{n}\to\dotsb\to v_{1}\xra{f_{1}}v}
\times (t_{n},\dotsc,t_{1}) \in \Gamma_{n+1}(u,v) \times I^{n},
$$
\noindent then 
%
%
\setcounter{equation}{\value{thm}}\stepcounter{subsection}
\begin{equation}\label{eone}
f_{1}\ccsub{t_{1}}f_{2}\ccsub{t_{n-1}}\dotsb\ccsub{t_{n-1}}
f_{n}\ccsub{t_{n}}f_{n+1}
\sim 
f_{1} \circ_{t_1} \dotsb \circ_{t_{i-1}} 
(f_i f_{i+1}) \circ_{t_{i+1}} \dotsb
\circ_{t_n} \circ f_{n+1}
\end{equation}
\setcounter{thm}{\value{equation}}
\noindent if \ $t_{i} =0$ \ for some \ $1 \leq i \leq n$. 

In \eqref{eone}, \ $(f_{i} f_{i+1})$ \ denotes \ $f_{i}$ \ composed with \ 
$f_{i+1}$. \ Loosely speaking, the relation $\sim$ means 
that \ $\circ_{t}$ \ becomes, for \ $t=0$, \ an ordinary composition of
morphisms in $\Gamma$.
The categorial composition in \ $\WG$ \ is given by the concatenation:
$$
(f_{1}\circ_{t_{1}} \dotsb \circ_{t_{l}} f_{l+1}) \circ 
(g_{1}\circ_{u_{1}} \dotsb \circ_{u_{k}} g_{k+1})\DEF 
(f_{1}\circ_{t_{1}} \dotsb \circ_{t_{l}} f_{l+1} \circ_{1} g_{1}
\circ_{u_{1}} \dotsb \circ_{u_{k}} g_{k+1}). 
$$
\end{defn}

Definition \ref{dwg} is based on the definition given in 
\cite[p.\ 367]{SchVo}; our version is simpler because we work without units.

\subsection{The cubical structure}
\label{scs}\stepcounter{thm}

The category \ $\WG$ \ is enriched not only over topological spaces, but
also cubically \ -- \  with a cubical set \ $\Hom_{\WG}(u,v)$ \ 
associated to each pair of objects \ $u,v$ \ of $\Gamma$  as follows:

For each \ $(n+1)$-chain
$$
\Phi=\lra{u \xra{f_{n+1}}v_{n}\xra{f_{n}}\dotsb\xra{f_{2}}v_{1}\xra{f_{1}} v}
$$
\noindent from $u$ to $v$, we have an $n$-dimensional cube \
$C[\Phi]$ \ of equivalence classes of points of \ 
$\{\Phi\} \times I^{n}\subset \Gamma_{n+1}(u,v) \times I^{n}$. \ Clearly,

\begin{enumerate}
\renewcommand{\labelenumi}{(\alph{enumi})}
\item the $k$-th $0$-facet of \ $C[\Phi]$ \ ($1\leq k\leq n$) \ equals to \ 
$C[\Phi^{(k)}]$, \ where \ $\Phi^{(k)}$ \ is obtained from $\Phi$ by composing 
the maps at the $k$-th internal node of $\Phi$;
\item the $k$-th $1$-facet of \ $C[\Phi]$ \ equals to the product of the 
two cubes \ $C[\Phi'_{k}]$ \ and \ $C[\Phi''_{k+1}]$ \ of dimensions \ 
$k-1$ \ and \ $n-k$ \ respectively, where \ 
$$
\Phi'_{k}\DEF \lra{v_{k}\xra{f_{k}}v_{k-1}\dotsb v_{1}\xra{f_{1}} v}
$$ 
and
$$
\Phi''_{k}\DEF\lra{u\xra{f_{n+1}}v_{n}\dotsb v_{k+1}\xra{f_{k+1}}v_{k}}.
$$
\end{enumerate}

Each facet of the cube \ $C[\Phi]$ \ has the form
$\{f_{1}\csub{t_{1}}f_{2}\dotsb\csub{t_{n}}f_{n+1}~|\ t_{i} \in \{0,1\} \}$ \ 
for some fixed \ $1\leq i\leq n$, \ and a
$k$-dimensional face of \ $C[\Phi]$ \ is defined by 
requiring \ $(n-k)$ \ of the parameters \ $(t_{1},\dotsc,t_{n})$ \ to take a 
fixed value in \ $\{0,1\}$.

In particular, vertices of \ $C[\Phi]$ \ are points
for which all \ $t_{i}\in \{0,1\}$. \ Using relation \eqref{eone}, we can
remove from $f_{1}\csub{t_{1}}f_{2}\dotsb\csub{t_{n}}f_{n+1}$ 
all operations \ $\circ_{0}$ \ by composing the adjacent morphisms. 
Therefore the vertices of the cubes of \ $\WG$ \  are in fact indexed by 
chains of composable morphisms. For example, the vertex \ 
$f_{1}\circ_{1}f_{2}\circ_{0}f_{3}\circ_{1}f_{4}$ \ of \ $C[\Phi]$ \ with
$$
\Phi = \lra{ u \stackrel{f_4}{\longrightarrow}
v_{3} \stackrel{f_{3}}{\longrightarrow} v_{2} \stackrel{f_{2}}{\longrightarrow}
v_{1} \stackrel{f_{1}}{\longrightarrow} v }
$$
is indexed by the $3$-chain \ 
$f_{1}(f_{2}f_{3})f_{4}$, \ vertex \ 
$f_{1}\circ_{0}f_{2}\circ_{0}f_{3}\circ_{1}f_{4}$ \ by the
$2$-chain \ $(f_{1}f_{2}f_{3})f_{4}$, \ and so on.

The (cubical) \emph{$k$-skeleton} of \ $\WG$, \ denoted by \ $\sk{k}\WG$, \ is 
the subcategory generated by the union of all faces of dimension $\leq k$. \ 

There is an obvious augmentation functor \ $\varepsilon:\WG\to\Gamma$, \ 
with a ``section'' \ $\sigma:\Gamma\to\WG$, \ which actually lands in
the $0$-skeleton \ $\sk{0}\WG$. \ The fiber \ $\varepsilon^{-1}(f)$ \ 
is contractible for each map $f$ of $\Gamma$.

\begin{remark}\label{runic}\stepcounter{subsection}
We will be interested mainly in the single cubical 
set \ $\Hom_{\WG}(\vi,\vf)$, \ which will be called the 
\emph{total mapping space} of \ $\WG$, \ and denoted by \ $\wG$, \ or 
simply $\PP$. \ However, we do need the full 
structure of \ $\WG$ \ as a cubically enriched category \ -- \ and this is
expressed by relations among the various subcomplexes of \ $\PP$.
\end{remark}

Let us abbreviate the notation of Definition \ref{dwg} by writing \ 
$fg$ \ for the composition \ $f\circ_{0}g$, \ $(f)(g)$ \ for the
composable sequence \ $f\circ_{1}g$, \ and \ $f\circ g$ \ for the
$1$-cube \ $f\circ_{t}g$ \ ($0\leq t\leq 1$) \ -- \ and similarly in
higher dimensions.

\begin{example}\label{esqur}\stepcounter{subsection}  
For the lattice \ $L_{3}$ \ (\S~\ref{elat}) with three composable morphisms 
$\vi \xra{h} v_3 \xra{g}v_2 \xra{f}\vf$, \ we see that \
$WL_{3}(\vi,\vf)$ \ is the square depicted in Figure~\ref{fig1},
where the arrows, indicating the direction \ $0\to 1$, \ are included
to identify more easily the two faces of each sub-cube.

\setcounter{figure}{\value{thm}}\stepcounter{subsection}
\begin{figure}[htbp]
%
%
\begin{center}
\begin{picture}(150,80)(0,0)
\put(50,75){\circle*{5}}
\put(5,71){$(fg)(h)$}
\put(115,75){\vector(-1,0){60}}
\put(70,81){{\scriptsize{$(f\circ g)(h)$}}}
\put(50,70){\vector(0,-1){60}}
\put(13,40){{\scriptsize{$(fg)\circ h$}}}
\put(120,75){\circle*{5}}
\put(125,71){$(f)(g)(h)$}
\put(120,70){\vector(0,-1){60}}
\put(123,40){{\scriptsize{$(f)(g\circ h)$}}}
\put(120,5){\circle*{5}}
\put(125,1){$(f)(gh)$}
\put(115,5){\vector(-1,0){60}}
\put(70,10){{\scriptsize{$f\circ(gh)$}}}
\put(50,5){\circle*{5}}
\put(15,2){$(fgh)$}
\put(67,40){\framebox{{\scriptsize{$f\circ g\circ h$}}}}
\end{picture}
\end{center}
\caption{\label{fig1}The square \ $WL_{3}(\vi,\vf)$}
\end{figure}
\setcounter{thm}{\value{figure}}
\end{example}

\begin{example}\label{ecube}\stepcounter{subsection}  
Similarly, for the lattice \ $L_{4}$ \ with four composable morphisms \ 
$\vi \xra{k} v_{4} \xra{h} v_{3} \xra{g} v_{2} \xra{f} \vf$, \ we find that \ 
$WL_{4}(\vi,\vf)$ \ is the cube in Figure~\ref{fig2}.

\setcounter{figure}{\value{thm}}\stepcounter{subsection}
\begin{figure}[htbp]
%
%
\begin{center}
\begin{picture}(380,200)(0,0)
%
%
\put(40,75){\circle*{3}}
\put(0,79){{\scriptsize{$(fgh)(k)$}}}
\put(135,75){\vector(-1,0){90}}
\put(65,66){{\scriptsize{$(fg\circ h)(k)$}}}
\put(40,70){\vector(0,-1){60}}
\put(1,40){{\scriptsize{$(fgh)\circ k$}}}
\put(140,75){\circle*{3}}
\put(102,80){{\scriptsize{$(fg)(h)(k)$}}}
\put(140,70){\vector(0,-1){60}}
\put(142,37){{\scriptsize{$(fg)(h\circ k)$}}}
\put(140,5){\circle*{3}}
\put(145,0){{\scriptsize{$(fg)(hk)$}}}
\put(135,5){\vector(-1,0){90}}
\put(75,-4){{\scriptsize{$fg\circ hk$}}}
\put(40,5){\circle*{3}}
\put(7,2){{\scriptsize{$(fghk)$}}}
\put(65,40){\framebox{{\scriptsize{$fg\circ h\circ k$}}}}
%
%
\put(225,175){\circle*{3}}
\put(180,175){{\scriptsize{$(f)(gh)(k)$}}}
\put(220,170){\vector(-2,-1){180}}
\put(80,125){{\scriptsize{$(f\circ gh)(k)$}}}
\put(325,175){\circle*{3}}
\put(330,176){{\scriptsize{$(f)(g)(h)(k)$}}}
\put(320,175){\vector(-1,0){90}}
\put(240,180){{\scriptsize{$(f)(g\circ h)(k)$}}}
\put(325,170){\vector(-2,-1){180}}
\put(245,122){{\scriptsize{$(f\circ g)(h)(k)$}}}
\put(139,116){\framebox{{\scriptsize{$(f\circ g\circ h)(k)$}}}}
\put(145,104){\scriptsize{top face}}
%
%
\multiput(225,170)(0,-3){23}{\circle*{.5}}
\put(225,105){\vector(0,-1){3}}
\put(180,135){{\scriptsize$(f)(gh\circ k)$}}
\put(325,165){\vector(0,-1){62}}
\put(328,135){{\scriptsize$(f)(g)(h\circ k)$}}
\put(325,100){\circle*{3}}
\put(328,93){{\scriptsize{$(f)(g)(hk)$}}}
\multiput(325,100)(-3,0){31}{\circle*{.5}}
\put(231,100){\vector(-1,0){3}}
\put(250,103){{\scriptsize$(f)(g\circ hk)$}}
\put(225,100){\circle*{3}}
\put(223,90){{\scriptsize$(f)(ghk)$}}
\put(228,160){\framebox{{\scriptsize{$(f)(g\circ h\circ k)$}}}}
%
%
\put(325,100){\vector(-2,-1){180}}
\put(262,60){{\scriptsize$(f\circ g)(hk)$}}
\multiput(225,100)(-2,-1){55}{\circle*{.5}}
\multiput(85,30)(-2,-1){20}{\circle*{.5}}
\put(45,10){\vector(-2,-1){3}}
\put(191,77){{\scriptsize$f\circ ghk$}}
\put(170,55){\framebox{{\scriptsize{$f\circ g\circ hk$}}}}
%
%
\put(35,163){\scriptsize{left face}}
\put(30,150){\framebox{{\scriptsize{$f\circ gh\circ k$}}}}
\put(50,145){\line(0,-1){20}}
\put(50,125){\vector(1,-1){15}}
%
%
\put(290,17){\scriptsize{right face}}
\put(280,30){\framebox{{\scriptsize{$(f\circ g)(h\circ k)$}}}}
\put(330,40){\line(0,1){20}}
\put(330,60){\vector(-1,1){20}}
\end{picture}
\end{center}
\caption{\label{fig2}The cube \ $\PP = WL_{4}(\vi,\vf)$}
\end{figure}
\setcounter{thm}{\value{figure}}
\end{example}

\begin{defn}\label{dbas}\stepcounter{subsection}  
The category \ $\BG$ \ is defined to be the cubical subcategory of \ 
$\WG$ \ with \ $\Obj \BG=\Obj\WG=\Obj\Gamma$, \ and with morphisms 
given by \ $\BG(u,v)\DEF\WG(u,v)$ \ if \ $(u,v)\neq(\vi,\vf)$, \ while
$$
\BG(\vi,\vf):=
\bigcup \left\{\alpha\circ\beta\ | \ \beta\in\WG(\vi,w),\
\alpha\in\WG(w,\vf),\ \vi\neq w\neq\vf
\right\}.
$$
\noindent As usual, we  abbreviate \ $\Hom_{\BG}(u,v)$ \ to \ $\BG(u,v)$ \ 
and  \  \ $\Hom_{\WG}(u,v)$ \ to \ $\WG(u,v)$.
\end{defn}

Thus \ $\BG(\vi,\vf)$ \ consists of all decomposable morphisms, \ 
so that as a cubical set its facets are of the form \ 
$C'_{k}[\Phi]\DEF\{(t_{1},\dotsc,t_{n})\in C[\Phi]~:\ t_{k}=1\}$ \ for
some fixed \ $1 \leq k \leq n$ \ and maximal chain $\Phi$. \ We will
sometimes call these the \emph{basic facets} of \ $\WG(\vi,\vf)$, \
and \ $\BG(\vi,\vf)$ \ will be called the \emph{basis} of \ $\WG$. \
We denote   $\BG(\vi,\vf)$ \ by \ $b\PP$ \
and consider it as a (cubical) subspace of $\PP$. 
For instance, in Example \ref{esqur} above, the basis of \ $WL_{3}$ \ consists
of the top and right edges of the square in Figure \ref{fig1}.

\begin{remark}\label{rzero}\stepcounter{subsection}
The $0$-skeleton \ $\sk{0}\WG$ \ is the free nonunital category on the set
of arrows of $\Gamma$.  Therefore each arrow \ $f \in \Gamma$ \ is
replicated in \ $\WG$, \ but the inclusion \ $\sigma : \Gamma \to \WG$ \ 
defined in this way is \emph{not} a functor. 
So for \ $u,v \in V$, \ we will sometimes identify elements $f \in 
\Gamma(u,v)$ with their images $\sigma(f) \in W\Gamma(u,v)$, that is,  
consider \ $\Gamma(u,v)$ \ as a 
\emph{subset} of \ $\WG(u,v)$, \ and emphasize this by writing \ 
$\Gamma \subset \sk{0}\WG$.

Similarly, \ $\sk{0}\BG$ \ is the free nonunital category on the set
of arrows of \ $\Gamma \setminus \{\fimax\}$. 
\end{remark}

%
%
\begin{prop}\label{pone}\stepcounter{subsection}  
For any lattice $\Gamma$, \ $\PP=\WG(\vi,\vf)$ \ is 
combinatorially isomorphic to the cone \ $Cb\PP$ \ on its basis \ 
$b\PP := \BG(\vi,\vf)$, \ with vertex of the cone corresponding to the unique
maximal $1$-chain \ $\lra{\vi\xra{\fimax}\vf}$.
\end{prop}

\begin{proof}
Recall that \ $C\bPP$ \ is the quotient \ 
$(I \times b\PP)/(\{0\} \times b\PP)$. \ 
Define a continuous map \ $\alpha : Cb\PP \to \PP$ \ as follows. 

Let \ $s \times [f_1 \circ_{t_1} f_2 \dotsb f_n \circ_{t_n} f_{n+1}]$ \ 
be a point of \ $I \times b\PP$, \ where 
$$
\langle\vi \stackrel {f_{n+1}}{\vlra}
v_n \stackrel {f_{n}}{\longrightarrow} \dotsb \stackrel 
{f_{2}}{\longrightarrow} v_1 \stackrel {f_{1}}{\longrightarrow}\vf\rangle
$$
is a maximal chain, \ $s,t_{1},\dotsc,t_{n} \in I$, \  
and \ $[-]$ \ denotes the equivalence class, as usual.
Let 
$$
\tilde \alpha(s \times [f_1 \circ_{t_1} f_2 \dotsb f_n \circ_{t_n}
f_{n+1}])
:= [f_1 \circ_{st_1} f_2 \dotsb f_n \circ_{st_n} f_{n+1}] \in \PP.
$$
\noindent If \ $s=0$, \ clearly 
$$
\tilde \alpha(s \times [f_1 \circ_{t_1} f_2 \dotsb f_n \circ_{t_n}
f_{n+1}]) = [f_1 \circ_{0} f_2 \dotsb f_n \circ_{0} f_{n+1}] = [\fimax].
$$
\noindent Thus \ 
$\tilde \alpha (\{0\}\times b\PP) = \{\fimax\}$, \ so that \ 
$\tilde\alpha$ \ induces a map \ $\alpha : Cb\PP \to \PP$ \ which sends
the vertex of the cone to \ $[\fimax] \in \PP$.

The inverse of $\alpha$ can be constructed as follows: let \ 
$f_{1}\circ_{t_1} f_{2} \dotsb f_{n} \circ_{t_n} f_{n+1}$ \ represent 
a point of $\PP$, and let \ $s := \max\{t_{1},\dotsc,t_{n}\}$. \ If \ 
$s=0$ \ (which means that \ $t_{i}=0$ \ for all $i$), let  \ 
$\beta([f_{1} \circ_{t_1} f_{2} \dotsb f_{n} \circ_{t_n}f_{n+1}])$ \ 
equal to the vertex of \ $Cb\PP$. \ If \ $s > 0$, \ then 
$$
\beta(f_{1}\circ_{t_1}f_{2} \dotsb f_{n} \circ_{t_n}
f_{n+1}) : = \pi(s \times [f_1 \circ_{t_1/s} f_2 \dotsb f_n \circ_{t_n/s}
f_{n+1}]) \in Cb\PP,
$$
\noindent where \ $\pi : I \times b\PP \to Cb\PP$ \ is the projection.
Observe that at least one of \ $t_1/s,\dotsc,t_n/s$ \ equals $1$, 
so indeed \ $[f_1 \circ_{t_1/s} f_2 \dotsb f_n \circ_{t_n/s}f_{n+1}] \in b\PP$.

It is easily verified that the map \ $\beta:\PP \to Cb\PP$ \ is
well-defined, and inverse to $\alpha$. 
\end{proof}

\begin{example}\label{ebasis}{\stepcounter{subsection}}
In Example \ref{ecube} above, the basis of \ $WL_{4}$ \ is the union
of the three $2$-dimensional faces opposite the vertex \ $(fghk)$, \ as
depicted in Figure \ref{fig3}; and the cube of Figure \ref{fig2} is indeed 
homeomorphic to the cone on \ $bWL_{4}$.

\setcounter{figure}{\value{thm}}\stepcounter{subsection}
\begin{figure}[htbp]
%
%
\begin{center}
\begin{picture}(260,160)(0,0)
%
%
\put(140,140){\circle*{3}}
\put(120,147){{\scriptsize{$(f)(gh)(k)$}}}
\put(136,140){\vector(-1,0){72}}
\put(73,132){{\scriptsize{$(f\circ gh)(k)$}}}
\put(60,140){\circle*{3}}
\put(22,142){{\scriptsize{$(fgh)(k)$}}}
\put(60,137){\vector(0,-1){134}}
\put(15,70){{\scriptsize{$(fg\circ h)(k)$}}}
\put(60,0){\circle*{3}}
\put(13,0){{\scriptsize{$(fg)(h)(k)$}}}
\put(140,75){\circle*{3}}
\put(87,75){{\scriptsize{$(f)(g)(h)(k)$}}}
\put(140,80){\vector(0,1){55}}
\put(88,115){{\scriptsize{$(f)(g\circ h)(k)$}}}
\put(136,72){\vector(-1,-1){69}}
\put(97,25){{\scriptsize{$(f\circ g)(h)(k)$}}}
\put(70,95){\framebox{{\scriptsize{$(f\circ g\circ h)(k)$}}}}
%
%
\put(145,140){\vector(1,0){70}}
\put(155,132){{\scriptsize{$(f)(gh\circ k)$}}}
\put(220,140){\circle*{3}}
\put(225,141){{\scriptsize{$(f)(ghk)$}}}
\put(220,75){\circle*{3}}
\put(225,70){{\scriptsize{$(f)(g)(hk)$}}}
\put(220,80){\vector(0,1){55}}
\put(222,110){{\scriptsize{$(f)(g\circ hk)$}}}
\put(145,75){\vector(1,0){70}}
\put(155,80){{\scriptsize{$(f)(g)(h\circ k)$}}}
\put(152,105){\framebox{{\scriptsize{$(f)(g \circ h\circ k)$}}}}
%
%
\put(67,0){\vector(1,0){149}}
\put(135,5){{\scriptsize{$(fg)(h\circ k)$}}}
\put(220,0){\circle*{3}}
\put(225,0){{\scriptsize{$(fg)(hk)$}}}
\put(220,72){\vector(0,-1){68}}
\put(222,35){{\scriptsize{$(f\circ g)(hk)$}}}
\put(145,45){\framebox{{\scriptsize{$(f\circ g)(h\circ k)$}}}}
\end{picture}
\end{center}
\caption{\label{fig3}Basis for\ $WL_\protect{4}(\vi,\vf)$}
\end{figure}
\setcounter{thm}{\value{figure}}
\end{example}

We now list some properties of the $W$-construction \ $\WG$ \ which
will be needed in Section \ref{chho}, using the following terminology:

For a finite set $V$ (which will be the set of objects of the category 
$\Gamma$), a \emph{$V$-graded space} will be a sequence \ 
$\X = \{X(v)\ | \ v \in V\}$ \ of  well-pointed topological spaces; 
a \emph{homotopy equivalence} between $V$-graded spaces $\X$ and $\Y$
is a sequence \ $\bh= \{ h_{v} : X(v)\to Y(v)~|\ v\in V\}$ \ of homotopy
equivalences in \ $\Ta$; \ if such an $\bh$ exists,  we say $\X$ and
$\Y$ have the same \emph{homotopy type}.

If $\Xi$ is a (topological) category with \ $\Obj(\Xi) = V$, \ 
a $V$-graded space $\X$ is called a \emph{$\Xi$-space} if there is a
\emph{$\Xi$-structure} on $\X$ \ -- \ that is, a (continuous) functor \ 
$G :\Xi \to\Ta$ \ such that \ $G(v) = X(v)$ \ for all \ $v \in V$.

The following statement, which we will need for the proof of 
Theorem~\ref{tone}, is a simplified version of \cite[Theorem~IV.4.37]{BVogHI} 
which was formulated there for an arbitrary colored operad. The proof
makes use of units, but applies by a slight modification in our case too.
%
%
\begin{thm}\label{tzero}\stepcounter{subsection}
A $V$-graded space $\X$ admits a \ $\WG$-structure if and only if it
is homotopy equivalent to a $\Gamma$-space. More precisely,

\begin{enumerate}
\renewcommand{\labelenumi}{(\alph{enumi})\ }
\item suppose that $\X$ is a \ $\WG$-space with the structure given by
a functor \ $\B:\WG\to\Ta$. \ Then there exists a $\Gamma$-space $\Y$ 
given by a functor \ $F : \Gamma \to\Ta$, \ and a homotopy equivalence \ 
$\bh: \X \xra{\sim}\Y$ \ such that the diagram
%
%
\setcounter{equation}{\value{thm}}\stepcounter{subsection}
\begin{equation}\label{efour}
\unitlength=1.2pt
\begin{picture}(100,30)(0,30)
\put(100,30){\makebox(0,0){\scriptsize $F(f)$}}
\put(-5,30){\makebox(0,0){\scriptsize $\B(\sigma(f))$}}
\put(50,3){\makebox(0,0){\scriptsize $h_v$}}
\put(50,57){\makebox(0,0){\scriptsize $h_u$}}
\put(90,10){\makebox(0,0){$F(v)$}}
\put(90,50){\makebox(0,0){$F(u)$}}
\put(10,10){\makebox(0,0){$\B(v)$}}
\put(10,50){\makebox(0,0){$\B(u)$}}
\put(90,40){\vector(0,-1){20}}
\put(10,40){\vector(0,-1){20}}
\put(22,10){\vector(1,0){56}}
\put(22,50){\vector(1,0){56}}
\end{picture}
\end{equation}
\setcounter{thm}{\value{equation}}
\vglue 25pt

\noindent is homotopy commutative for all arrows \ $f :u\to v$ \ 
of \ $\Gamma$.
\item Conversely, if $\Y$ is a $\Gamma$-space, with structure
given by \ $F: \Gamma \to\Ta$, \ and \ $\bh: \X \to \Y$ \ is a 
homotopy equivalence of $V$-graded spaces, then there is a \ 
$\WG$-structure \ $\B : \WG \to \Ta$ \ on $\X$ such that the 
diagram \eqref{efour} is homotopy commutative for all \ $f :u\to v$ \
as above.
\end{enumerate}
\end{thm}

The strict $\Gamma$-space $\Y$ whose existence is claimed in the first part
of the Theorem is called a \emph{rectification} of the \ 
$\WG$-space $\X$. The original result in \cite[Thm.~IV.4.37]{BVogHI}
in fact states that the homotopy equivalence $\bh$ can be
equipped with a structure of a strongly homotopy $\Gamma$-map,
(in an appropriate sense). The homotopy commutativity of~(\ref{efour}) then
follows from the existence of this structure. 

For a subcategory $U$ of \ $\WG$, \ $u,v \in V$, \ and \ $n \geq 0$, \ let \ 
$U_{n+1}(u,v) \subset \Gamma_{n+1}(u,v) \times I^{n}$ \ denote the
set of all elements representing morphisms in $U$.
Following~\cite{SchVo}, we call $U$ {\em admissible\/} if each morphism
of $U$ that decomposes in $\WG$, decomposes also in $U$, and if,
moreover, the inclusion 
$$
U_{n+1}(u,v) \cup \Gamma_{n+1}(u,v) \times \partial I^n
\hookrightarrow 
\Gamma_{n+1}(u,v) \times I^n
$$
is a closed cofibration for all \ $n,u,v$. 
Observe that both \ $\sk{0}\WG$ \ and \ $\sk{0}b\WG$ \ are admissible
subcategories of \ $\WG$; \ 
in both cases \ $U_{n+1}(u,v)\subset\Gamma_{n+1}(u,v)\times\partial I^{n}$.

The above condition is a bit weaker than the one in~\cite[p.~374]{SchVo}, 
because we do not have units. The following statement is a version
of~\cite[Prop.~4.4]{SchVo}, which will be used in the proof of
Proposition~\ref{pfour}:

%
%
\begin{prop}\label{ptwo}\stepcounter{subsection}  
Let $U$ be an admissible subcategory of \ $\WG$. \ Suppose we are
given a functor \ $F_{0}: \WG \to \Ta$ \ and a homotopy through
functors \ $H_{t} : U \to \Ta$ \ such that \ $H_{0} = F_{0}\rest{U}$. \ 
Then there exists an \ $F_{t}$ \ extending both \ $F_{0}$ \ and \ $H_{t}$.
\end{prop}

\begin{remark}\label{rsv}\stepcounter{subsection}  
There is a version of \ $\WG$, \ based on a construction of Cordier and Porter 
(in \cite[\S 2]{CPorV}), which is enriched over \emph{simplicial}, rather
than cubical, sets:

Given a lattice $\Gamma$ as above, the category \ $\WsG$ \ has the
same objects as $\Gamma$, and for each pair of nodes \ $(u,v)$ \ of $\Gamma$,
$\WsG(u,v)\in\Ss$ \ is a simplicial set with one $r$-simplex \ 
$\sigma_{(\U,\Phi)}$ \ for each chain \ 
$\Phi=\lra{u=v_{n+1}\xra{f_{n+1}}v_{n}\xra{f_{n}}\dotsb\xra{f_{2}}
v_{1}\xra{f_{1}}v_{0}=v}=f_{1}f_{2}\dotsb f_{n}f_{n+1}$ \ 
and each partition \ $\U=(U_{1},\dotsc,U_{r})$ \ of a subset \ 
$U\subset \{v_{1},\dotsc,v_{n}\}$ \ of the set of internal nodes of $\Phi$, \ 
with each of the sets \ $U_i$ \ nonempty. The faces of $\sigma$ are defined by:
\begin{enumerate}
\renewcommand{\labelenumi}{(\roman{enumi})}
\item $d_{0}(\sigma)\DEF\sigma_{(\U^{(0)},\Phi')}$, \ where \ $\Phi'$ \
is obtained from $\Phi$ \ by carrying out the compositions at each node \ 
$v_{i}\in U_{1}$,  \ and \ $\U^{(0)}\DEF(U_{2},\dotsc, U_{r})$.
\item $d_{j}(\sigma)\DEF\sigma_{(\U^{(j)},\Phi)}$, \ where \ 
$\U^{(j)}\DEF (U_{1},\dotsc,U_{j-1},U_{j}\cup U_{j+1},U_{j+2},\dotsc,U_{r})$ \ 
for \ $0<j<r$.
\item $d_{r}(\sigma)\DEF\sigma_{(\U^{(r)},\Phi)}$, \ where \ 
$\U^{(r)}\DEF(U_{1},\dotsc,U_{r-1})$.
\end{enumerate}

The degenerate simplices are obtained by allowing partitions with \ 
$U_{j}=\emptyset$, \ and the simplicial composition map is defined by 
concatentation of chains in the obvious way. 
Note that \ $\WsG(u,v)$ \ (or the corresponding simplicial complex) 
clearly provides a canonical triangulation of the cubical set \ $\WsG(u,v)$. 

The construction \ $\WsG$ \  can be thought of as a functorial bar 
resolution: \ 
$$
B(\sk{0}W,\sk{0}W_{0},\Gamma)
$$
\noindent (cf.\ \cite[\S 9]{MayG}).
\end{remark}
 
%
%
\sect{Higher homotopy operations}
\label{chho}

With the constructions of Section~\ref{clwc} at hand, we are now in a 
position to define our basic object of interest:

\begin{defn}\label{did}\stepcounter{subsection}  
\emph{Initial data} for a higher homotopy operation is a lattice $\Gamma$, 
together with a $\Gamma$-diagram up-to-homotopy \ -- \ 
i.e.,  a functor \ $\A:\Gamma\to\hoT$. 

Note that if \ $\sigma:\Gamma\to\WG$ \ is the (non-functorial) canonical
section of the augmentation \ $\varepsilon : \WG\to\Gamma$, then for every
continuous functor \ $F:\WG\to\Ta$, \ the composite \ 
$\pi\circ F\circ\sigma:\Gamma\to\hoT$ \ is a functor, and \ 
$\pi\circ F$ \ factors through $\varepsilon$, \ where \
$\pi:\Ta\to\hoT$ \ is the obvious projection functor. 

A \emph{rectification} of the initial data \ $\A:\Gamma\to\hoT$ \ 
is then a strict $\Gamma$-diagram realizing $\A$ \ -- \ i.e., a functor \ 
$F:\Gamma\to\Ta$ \ such that \ $\pi\circ F$ \ is naturally isomorphic
to $\A$.
\end{defn}

Recall that the (right) \emph{half-smash} \ $X\rtimes K$ \ of topological
spaces $X$ and $K$, where $X$ is pointed, is 
defined to be \ $(X\times K)/(\{\ast\}\times K)=X\wedge K_{+}$, \
where \ $K_{+}$ \ is $K$ with a disjoint basepoint added. \ 
$X\rtimes K$ \ is again a pointed space, with the class of \ 
$\{*\} \times K$ \ 
as the distinguished point. 

Given two pointed spaces \ $X,Y\in\Ta$ \ and any space \ $K\in\TT$, \ 
the following adjointness relation holds:
%
%
\setcounter{equation}{\value{thm}}\stepcounter{subsection}
\begin{equation}\label{etwo}
\Ta(X\rtimes K,Y) \cong \Ta(X\wedge K_{+},Y) \cong
\Ta( K_{+},\Ta(X,Y)) \cong \TT(K, \Ta(X,Y)),
\end{equation}
\setcounter{thm}{\value{equation}}

\noindent provided one works in the category of (not necessarily
connected) pointed compactly-generated spaces.

\begin{con}\label{cdist}\stepcounter{subsection}  
In \eqref{etwo}, we will typically have \ $X = \A(\vi)$, \ 
$Y =\A(\vf)$ \ and $K = b\PP$. \ In this case, we will make 
\emph{no distinction} between maps \ 
$\A(\vi)\rtimes b\PP \to \A(\vf)$ \ and maps\ $b\PP
\to\Ta(\A(\vi),\A(\vf))$.
\end{con}

\begin{defn}\label{dhho}\stepcounter{subsection}  
Given initial data \ $\A:\Gamma\to\hoT$, \ \emph{complete data} for
the corresponding higher homotopy operation consists of a continuous
functor \ $\BA:\BG\to\Ta$ \ 
such that \ $\pi\circ\BA  = \A\circ (\varepsilon\rest{\BG})$.

The corresponding \emph{higher order homotopy operation} \ is the subset \ 
$$
\llrr{\A}\subset [\A (\vi) \rtimes b\PP,\A (\vf)]_{\hoT}
$$
\noindent consisting of the homotopy equivalence classes of maps 
$$
\BA|_{b\WG(\vi,\vf)}: b\WG(\vi,\vf) = 
\bPP \longrightarrow\Ta(\A(\vi),\A(\vf))
$$
induced by all possible complete data \ $\BA$ \ for $\A$.
\end{defn}

\begin{defn}\label{dvan}\stepcounter{subsection}  
The higher operation \ $\llrr{\A}$ \ is said to \emph{vanish} if it 
contains the homotopy class of a constant map $\bPP
\longrightarrow\Ta(\A(\vi),\A(\vf))$. \
\end{defn}

The proof of the following statement is based on Proposition~\ref{pone}.

%
%
\begin{prop}\label{pthree}\stepcounter{subsection}  
The operation \ $\llrr{\A}$ \ vanishes if and only if there exists  
a continuous functor \ $\B : \WG \to\Ta$ \ such that 
%
%
\setcounter{equation}{\value{thm}}\stepcounter{subsection}
\begin{equation}\label{eeight}
\pi\circ \B = \A\circ\varepsilon.
\end{equation}
\setcounter{thm}{\value{equation}}
\end{prop}

\begin{proof}

The vanishing of \ $\llrr{\A}$ \ means, by definition, that there are complete
data \ $\BA : b\WG\to\Ta$ \ such that \ $\BA\rest{b\WG(\vi,\vf)}$ \ is
homotopic to a constant map. \ We then 
define \ $\B : \WG \to \Ta$ \ on objects by \ $\B(v) := \A(v)$, for each \ 
$v\in V$, \ and on morphisms as follows. 

Recall that, for \ $(u,v)\neq (\vi,\vf)$ \ we have \ 
$\WG(u,v)=b\WG(u,v)$, \ so in this case we simply put \ 
$\B\rest{\WG(u,v)}\DEF\B\rest{b\WG(u,v)}$. \ 
{}For \ $(u,v) = (\vi,\vf)$, \ let \ $\B\rest{\WG(\vi,\vf)}$ \ be an
arbitrary extension of \ $\BA\rest{b\WG(\vi,\vf)}$.
Such an extension exists,
because $W\Gamma(\vi,\vf)$ is, by Proposition~\ref{pone}, a cone over 
$bW\Gamma(\vi,\vf)$ and $\BA\rest{b\WG(\vi,\vf)}$ is, by assumption, 
homotopic to a constant map.
It is obvious that $\B$ defined in this way is a functor \ -- \ 
since all decomposable morphisms of \ $\WG$ \ belong to \ $b\WG$, \ 
there are no categorial constraints on the extension
$\B$. Equation~\eqref{eeight} is clearly also satisfied.

The opposite implication \ -- \ that the existence of a functor \ 
$\B:\WG\to\Ta$ satisfying~\eqref{eeight} implies the vanishing
of \ $\llrr{\A}$ \
 -- \ also follows directly from Proposition~\ref{pone}. 
\end{proof}

It is easy to see that~\eqref{eeight} is equivalent to
\[
\pi \circ \B(\sigma(f)) = \A(f),
\]
for each $f \in \Gamma$.
We can now formulate and prove the main result of this section:
%
%
\begin{thm}\label{tone}\stepcounter{subsection}
The homotopy operation \ $\llrr{\A}$ \ vanishes (and in particular, is 
defined) if and only if $\A$ has a rectification, so it is precisely 
the last obstruction to rectifying $\A$.
\end{thm}

\begin{proof}
The proof is an immediate consequence of Theorem~\ref{tzero}. 
Let $\X$ denote, throughout this proof, the $V$-graded space defined
by \ $X(v) := \A(v)$, \ for \ $v \in V$. 

Let us prove first that the
vanishing of \ $\llrr{\A}$ \ implies the existence of a rectification. By
Proposition~\ref{pthree}, vanishing of $\llrr{\A}$ 
is equivalent to the existence of a 
functor \ $\B : \WG \to \Ta$ \ such that \ $[\B(\sigma(f))] = \A(f)$ \ 
for each arrow \ $f: u \to v \in \Gamma$. \ 
The functor $\B$ is a \ $\WG$-structure on the $V$-graded space $\X$
and Theorem~\ref{tzero}(a) provides a $\Gamma$-space $\Y$ defined by a
functor \ $F : \Gamma \to \Ta$, \ together with a homotopy equivalence \ 
$\bh  = \{h_u : \B(u) \to F(u)~| \ u \in V\}$, \  
for which diagram~(\ref{efour}) commutes up to homotopy, that is
$$
[F(f)] \circ [h_u] = [h_v] \circ [\B(\sigma(f))]
$$
in \ $\hoT$. \ Since \ $[\B(\sigma(f))] = \A(f)$, \ 
this means that \ $F : \Gamma \to \Ta$ \ is a rectification of $\A$.

To prove that the existence of a rectification implies the
vanishing of \ $\llrr{\A}$, \ let \ 
$F : \Gamma\to \Ta$ \ be such a rectification \ -- \ that is, a $V$-graded
$\Gamma$-space \ $\Y = \{F(v)~| \ v \in V\}$, \ together with a system \ 
$\{h_u : \A(u) \to F(u)~| \ u \in V \}$ \ such that 
%
%
\setcounter{equation}{\value{thm}}\stepcounter{subsection}
\begin{equation}\label{ethree}
[F(f)] \circ [h_u] = [h_v] \circ \A(f)
\end{equation}%
\setcounter{thm}{\value{equation}}%
for each \ $f : u \to v \in \Gamma$. \ 
By Theorem~\ref{tzero}(b), there exists a \ $\WG$-structure $\B$ on $\X$ such
that \ $[F(f)] \circ [h_u] = [h_v] \circ [\B(\sigma(f))]$. \ 
Combining this with \eqref{ethree} and using the invertibility
of \ $[h_v]$ \ in \ $\hoT$ \ we see that \ $[\B(\sigma(f))] = \A(f)$ \ 
for each $f$, therefore \ $\llrr{\A} = 0$, \ by Proposition~\ref{pthree}.  
\end{proof}

\subsection{A pointed version}
\label{spv}\stepcounter{thm}

Since the category \ $\Ta$ \ is pointed, the null map $\ast$ has a 
special status. We may therefore require that whenever \ $\A(f)=\ast$ \ 
(the null class), \ we take \ $B{\A}(f)=\ast$ \ to be the null map itself
(rather than just nullhomotopic), and more generally that for any 
cubical face \ $C[\Psi]$ \ of \ $\WG$, \ if the corresponding map 
on the boundary \ $\partial C[\Psi]\times X\to Y$ \ is the null map, 
then the extension to \ $C[\Psi]\times X\to Y$ \ also be the null map. 
In this case we say that the corresponding complete data are \emph{pointed}.

This is in fact the most common sense in which the term ``secondary
(or higher order) homotopy operation'' is used. As we shall see in \S
\ref{rsmash} below, there are practical reasons why this is the 
prefered version.

Note that in this case one can actually replace the basis \ $\PP=\BG$ \ by a 
\emph{simplified basis} \ $\hat{\PP}=\hBG$, \ in which we collapse all
faces indexed by at least one zero map to a single point.  Of course,
the resulting cell complex may no longer be a polyhedron; see Figure 
\ref{fig4} in Example \ref{eltb} below.

\subsection{Relative operations}
\label{sro}\stepcounter{thm}

More generally, we may choose to specify a partial rectification of the 
diagram \ $\A:\Gamma\to\hoT$, \ in the sense that there is a sub-diagram \ 
$\Gamma'\subset\Gamma$ \ equipped with a specific lift of \ 
$\A\rest{\Gamma'}$ \ to a functor \ $F':\Gamma'\to\Ta$. \ 
In this case again we assume that the extensions to the appropriate faces 
of \ $b\PG$ \ are constant, and call the corresponding subset of \ 
$[\A(\vi)\rtimes b\PP,\A(\vf)]\text{~(rel~}\A(\vi)\rtimes b\PP_{\Gamma'}$) \ 
a \emph{relative} higher homotopy operation.

\begin{example}\label{etb}\stepcounter{subsection}  
The classical example of the Toda bracket fits into the variant of 
\S \ref{spv}, since in the usual description we require certain of the 
maps constituting the full data for this operation to be specifically the 
null map (rather than allowing any null-homotopic map, as one would in a 
fully homotopy-invariant description):

Recall that the usual Toda bracket is defined whenever one has three composable
maps \ $X \xra{\gamma} Y \xra{\beta} Z \xra{\alpha} W$, \ and each of the 
compositions \ $\alpha\circ\beta$ \ and \ $\beta\circ\gamma$ \ is 
nullhomotopic. We thus have a lattice $\Gamma$ as in Figure \ref{fig4},
with four commuting triangles marked \ 
$\fbox{F}$, \ $\fbox{G}$, \ $\fbox{H}$ \ and \ $\fbox{K}$. \ Note
that this is not a planar graph: the two outer (angled) edges, both marked \ 
$f_{\ast}$, \ should be identified.

\setcounter{figure}{\value{thm}}\stepcounter{subsection}
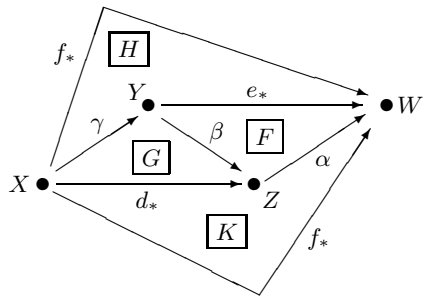
\begin{figure}[htb]
\begin{center}
%
%
\begin{picture}(150,120)(0,0)
%
%
\put(10,40){$\bullet$}
\put(0,40){{\scriptsize{$X$}}}
\put(18,48){\vector(3,2){31}}
\put(31,64){{\scriptsize{$\gamma$}}}
\put(18,43){\vector(1,0){70}}
\put(48,34){{\scriptsize{$d_{\ast}$}}}
\put(47,50){\framebox{\scriptsize{$G$}}}
\put(50,70){$\bullet$}
\put(45,75){{\scriptsize{$Y$}}}
\put(58,68){\vector(3,-2){30}}
\put(76,60){{\scriptsize{$\beta$}}}
\put(90,40){$\bullet$}
\put(96,34){{\scriptsize{$Z$}}}
%
%
\put(58,73){\vector(1,0){76}}
\put(90,76){{\scriptsize{$e_{\ast}$}}}
\put(97,44){\vector(3,2){38}}
\put(116,50){{\scriptsize{$\alpha$}}}
\put(140,70){$\bullet$}
\put(147,70){{\scriptsize{$W$}}}
\put(90,58){\framebox{\scriptsize{$F$}}}
%
%
\put(16,50){\line(1,3){20}}
\put(36,110){\vector(3,-1){100}}
\put(38,91){\framebox{\scriptsize{$H$}}}
\put(17,90){{\scriptsize{$f_{\ast}$}}}
%
%
\put(17,40){\line(2,-1){78}}
\put(95,1){\vector(2,3){42}}
\put(113,20){{\scriptsize{$f_{\ast}$}}}
\put(75,22){\framebox{\scriptsize{$K$}}}
\end{picture}
\end{center}
\caption{\label{fig4}Lattice for Toda bracket}
\end{figure}
\setcounter{thm}{\value{figure}}

The homotopy classes of
$\alpha,\beta$ and $\gamma$ as above define initial data $\A$ 
as indicated by the diagram, 
with \ $f_{\ast}$, \ $e_{\ast}$ \ and \ $d_{\ast}$ \ being null maps.  
The associated polyhedron \ $\PP$ \ is a square, as in Figure~\ref{fig5},
where \ $\BG (X,W)$ \ consists of the upper and right edges of the square. 

\setcounter{figure}{\value{thm}}\stepcounter{subsection}
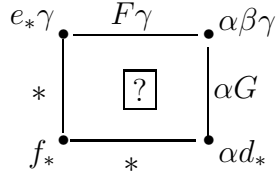
\begin{figure}[htbp]
\begin{center}
%
%
\begin{picture}(100,60)(0,0)
\put(20,10){{\scriptsize{$\bullet$}}}
\put(9,4){$f_{\ast}$}
\put(25,12){\line(1,0){47}}
\put(45,0){$\ast$}
\put(75,10){{\scriptsize{$\bullet$}}}
\put(81,4){$\alpha d_{\ast}$}
\put(22,15){\line(0,1){35}}
\put(10,27){$\ast$}
\put(20,50){{\scriptsize{$\bullet$}}}
\put(2,56){$e_{\ast}\gamma$}
\put(26,52){\line(1,0){47}}
\put(40,56){$F\gamma$}
\put(75,50){{\scriptsize{$\bullet$}}}
\put(81,54){$\alpha\beta\gamma$}
\put(77,12){\line(0,1){35}}
\put(79,28){$\alpha G$}
\put(45,27){\framebox{?}}
\end{picture}
\end{center}
\caption{\label{fig5}Square for Toda bracket}
\end{figure}
\setcounter{thm}{\value{figure}}
\noindent Since \  $\BG(X,Z)$ \ is the interval with endpoints \
$\beta \circ \gamma$ \ and \ $d_{\ast}$, \ and \ $\BG(Y,W)$ \ is the
interval with endpoints \ $\alpha \circ \beta$ \  and \ $e_{\ast}$, \ 
complete data are given by a choice of homotopies $F$ and $G$ 
corresponding to commutative triangles \ $\framebox{G}$ \ and \ 
$\framebox{F}$; \ the homotopies corresponding to triangles \ 
$\framebox{H}$ \ and \ $\framebox{K}$ \ are assumed to be trivial. 
The fill-in marked \ $\framebox{?}$ \ exists 
(with trivial homotopies along the edges marked ``$\ast$'') if and only if the 
corresponding Toda bracket \ $\llrr{\A}\subset[X\rtimes S^1,W]$ \ vanishes.
\end{example}

\subsection{Conditions for existence}
\label{sce}\stepcounter{thm}
Given initial data \ $\A : \Gamma \to \hoT$, \ it is natural to try to
construct complete data \ $\BA: b\WG \to \Ta$ \ and a functor \ $\B :
\WG \to\Ta$ \ inductively, using the skeletal filtration of \ $b\WG$.   

\begin{defn}\label{dnr}\stepcounter{subsection}  
For \ $n \geq 0$, \ an $n$-\emph{realization} of the initial data $\A$ 
is a functor \ $\bA_{n}:\sk{n}\BG\to\Ta$ \ such that \ 
$[\bA_{n}(f)]=\A (f)$ in \ $\hoT$ \ for any arrow $f$ in \ 
$\Gamma \setminus \{\fimax\} \subset \sk{0}b\WG$.
\end{defn}

It follows from the freeness of the category \ $\sk{0}\WG$ \ 
(\S \ref{rzero}) that there always exists a $0$-realization \ $\bA_{0}$. \  
It sends any partially parenthesized sequence 
representing a vertex of  \ $\sk{0}\BG$  --  such as \ $fg(hk\ell)m$  --
to the corresponding chain in $\Ta$ \ -- e.g., \ 
$\bA_{0}f \bA_{0}g \bA_{0}(h\circ k\circ \ell) \bA_{0}m$. \ 
Observe also that the initial data $\A$ are uniquely determined by any
of $0$-realizations of $\A$. 

Using the above terminology, we can also say that complete data are
given by a functor \ $\BA : b\WG \to \Ta$ \ extending some $0$-realization
of $\A$. In the same vein, the vanishing of \ $\llrr{\A}$ \ is
equivalent to the existence of a functor \ $\B:\WG\to\Ta$ \ that
extends some $0$-realization \ $\bA_{0}$. \  
Indeed, for such a functor we have \ $[\B(f)] = [\bA_{0}(f)]= \A(f)$ \
whenever \ $f \in \Gamma \setminus \{\fimax\}$, \ while \ 
$[\B(\fimax)] = \A(\fimax)$ \ follows from the fact that \ 
$\WG(\vi,\vf)$ \ is the cone over \ $b\WG(\vi,\vf)$ \ with vertex \ 
$\lra{\fimax}$. \ Compare also Proposition~\ref{pthree} and its proof.

Let us show that the  success or failure of the induction does not
depend on the first step of the construction:
%
%
\begin{prop}\label{pfour}\stepcounter{subsection}  
Suppose that there exists a $0$-realization \ $\bA_{0}:\sk{0}b\WG \to \Ta$ \ 
which can be extended to a functor \ $\B : \WG \to \Ta$. \ Then an
arbitrary $0$-realization \ $\bA_{0}'$ \ can be extended to some functor \ 
$\B':\WG \to \Ta$.
\end{prop}

\begin{proof}
Let \ $\bA_{0}:\sk{0}\BG\to\Ta$ \ be a $0$-realization of $\A$, \ 
$\B:\WG\to\Ta$ \ its functorial extension
and \ $\bA_{0}':\sk{0}\BG \to \Ta$ \  another $0$-realization of
$\A$. Since, by definition, 
%
%
\setcounter{equation}{\value{thm}}\stepcounter{subsection}
\begin{equation}\label{efive}
[\bA_{0}(f)] = [\bA_{0}'(f)] = \A(f)
\end{equation}
\setcounter{thm}{\value{equation}}
\noindent in \ $\hoT$, \ for each \ 
$f \in\Gamma \setminus \{\fimax\} \subset \sk{0}b\WG$, \ 
there exists a homotopy through functors \ $H_t:\sk{0}b\WG\to\Ta$ \ 
from \ $\bA_{0}$ \ to \ $\bA_{0}'$. \ Indeed, choose \ $H_{t}$ \ 
arbitrarily on generators \ $f \in\Gamma \setminus \{\fimax\}$ \ -- \ 
this is possible by \eqref{efive} \ -- \ and then extend \ $H_t$ \ 
functorially using the freeness of \ $\sk{0}b\WG$ \ (Remark~\ref{rzero}).   

Because \ $\sk{0}b\WG$ \ is an admissible subcategory of $\WG$, we 
may apply Proposition~\ref{ptwo} to obtain an
extension \ $F_{t}:\WG\to\Ta$ \ of $\B$ and \ $H_{t}$. \ Then \ 
$\B':= F_1:\WG\to\Ta$ \ is an extension of \ $\bA'$. \ 
\end{proof}

\begin{remark}\label{rsmash}\stepcounter{subsection}
In the case of a pointed higher homotopy operation, in the sense of 
\S \ref{spv}, we are often in the situation of Example \ref{etb}, in 
that one (or more) of the vertices of \ $b\PP$ \ is taken by $\A$ to the null
map. This implies that we may in fact think of the higher homotopy operation \ 
$\llrr{\A}$ \ as comprising a subset of \ 
$[\A \vi \wedge b\PP,\A \vf]_{\hoT}$ \  -- \ where we choose the above
vertex of \ $b\PP$ \ as the base point. This is useful when \ 
$b\PP$ \ is a sphere, up to homotopy, since in that case the higher operations
take value in an (abelian) group.
\end{remark}

\begin{example}\label{eltb}\stepcounter{subsection}  
Longer sequences of composable maps yield higher Toda brackets, as in
\cite{GWalkL,MoriHT}. For instance, given the lattice \ 
$L_{4}=\lra{\vi \xra{k} v_{4} \xra{h} v_{3} \xra{g} v_{2} \xra{f} \vf}$ \ of 
Example \ref{ecube}, we obtain a three-dimensional cube as in 
Figure \ref{fig2} for \ $WL_{4}$, \ with basis as in Figure \ref{fig3}.
However, if we assume that \ $f\circ g$, \ $g\circ h$, \ 
and \ $h\circ k$ \ are all null, and we are interested in the pointed
third-order operation (\S \ref{spv}), we may replace \ 
$\PP=bWL_{4}$ \ by the simplified basis \ $\widehat{\PP}=\hBG$ \ of 
Figure \ref{fig9}.

\setcounter{figure}{\value{thm}}\stepcounter{subsection}
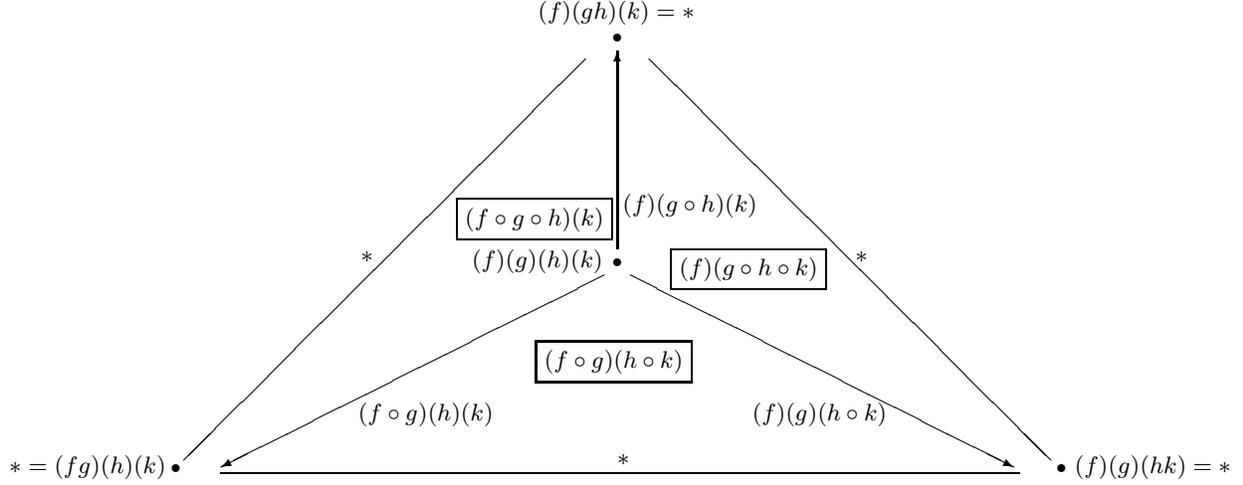
\begin{figure}[htbp]
%
%
\begin{center}
\begin{picture}(430,180)(0,0)
%
%
\put(230,165){\circle*{3}}
\put(200,172){{\scriptsize{$(f)(gh)(k)=\ast$}}}
\put(218,157){\line(-1,-1){152}}
\put(133,80){{\scriptsize{$\ast$}}}
\put(63,2){\circle*{3}}
\put(0,0){{\scriptsize{$\ast=(fg)(h)(k)$}}}
\put(230,80){\circle*{3}}
\put(175,78){{\scriptsize{$(f)(g)(h)(k)$}}}
\put(230,85){\vector(0,1){75}}
\put(232,99){{\scriptsize{$(f)(g\circ h)(k)$}}}
\put(225,75){\vector(-2,-1){145}}
\put(132,20){{\scriptsize{$(f\circ g)(h)(k)$}}}
\put(169,94){\framebox{{\scriptsize{$(f \circ g \circ h)(k)$}}}}
%
%
\put(242,157){\line(1,-1){152}}
\put(320,80){{\scriptsize{$\ast$}}}
\put(398,2){\circle*{3}}
\put(403,0){{\scriptsize{$(f)(g)(hk)=\ast$}}}
\put(235,75){\vector(2,-1){145}}
\put(281,20){{\scriptsize{$(f)(g)(h\circ k)$}}}
\put(250,75){\framebox{{\scriptsize{$(f)(g\circ h\circ k)$}}}}
%
%
\put(80,0){\line(1,0){302}}
\put(230,3){{\scriptsize{$\ast$}}}
\put(199,40){\framebox{{\scriptsize{$(f\circ g)(h\circ k)$}}}}
\end{picture}
\end{center}
\caption{\label{fig9}Simplified basis for \ $\widehat{bWL}_{4}$}
\end{figure}
\setcounter{thm}{\value{figure}}

If we carry out the indicated identifications (i.e., collapse the
boundary of the outer triangle to a point), we obtain a $2$-sphere, so we
are indeed in the situation describe in Remark \ref{rsmash}. 
\end{example} 

\subsection{Extending {$k$}-realizations}
\label{sekr}\stepcounter{thm}

Assume given initial data \ $\A:\Gamma\to\hoT$; then we always have
some extension of any $0$-realization \ $\bA_{0}:\sk{0}\BG\to\Ta$ \ to \ 
$\sk{1}\BG$, \ since this merely involves \emph{choosing} homotopies 
between maps which are homotopic by definition.
However, this is not true of higher skeleta:  a given \
$(k-1)$-realization may not extend to a $k$-realization, and it may
extend in more than one way; as expected, we have an obstruction theory.
Rather than describing it in detail, we shall sketch only the part
that is relevant to higher homotopy operations; in particular, we
disregard the obstructions for distinguishing between different extensions.

Note that, given what we know about \ $\bA_{0}$, \ it is natural to make
two requirements for our induction:

\begin{assume}\label{aekr}\stepcounter{subsection}
At the $k$-th stage of the induction \ ($k\geq 1$), \ we assume that:

\begin{enumerate}
\renewcommand{\labelenumi}{(\alph{enumi})}
\item we have fixed a \ $(k-1)$-realization \ $\bA_{k-1}:\sk{k-1}\BG\to\Ta$;
\item $\bA_{k-1}$ \ \emph{can} be extended over \ $\sk{k}\BG$ \  -- \ 
although we do \emph{not} choose a particular $k$-realization \ $\bA_{k}$.
\end{enumerate}
\end{assume}

Note also that if \ $n=\dim\BG$, \ the $n$-th induction assumption is
precisely what we want in order for \ $\llrr{\A}$ \ to be defined; 
the various choices of \ $\bA_{n}$ \ (for all possible choices of \ 
$\bA_{n-1}$) \ then define its members.

\begin{defn}\label{dac}\stepcounter{subsection}  
Given a maximal chain \ 
$$
\Phi=\lra{\vi =v_{n}\xra{f_{n}} v_{n-1}\xra{f_{n-1}}\dotsb
\xra{f_{2}} v_{1}\xra{f_{1}} v_{0}=\vf}
$$
\noindent of $\Gamma$, \ every basis point \ 
$\lra{f_{1}\csub{t_{1}}f_{2}\csub{t_{2}}\dotsb f_{n-1}\csub{t_{n-1}}f_{n}}$ \ 
in the corresponding subcomplex \ $C[\Phi]\cap \BG$ \ is in \ 
$\sk{k}\BG$ \ if and only if it has at least \ $(n-k-1)$ \ 
coordinates equal to $1$ or $0$. \ Each choice of a subset \ 
$M\subset\{1,2,\dotsc, n-1\}$ \ (those of the co\-ordinates which we 
wish to set equal to $1$) yields a partition of \ $\Phi$ \ into $r$ 
disjoint subchains: \ $\Phi=\Psi^{1}\Psi^{2}\dotsb\Psi^{r}$. \ 

To each proper subchain \ $\Psi^{i}=
\lra{ \vi^{i} \xra{f^{i}_{n_{i}+1}}v^{i}_{n_{i}}\xra{f^{i}_{n_{i}}}
\dotsb\xra{f^{i}_{2}}v^{i}_{1}\xra{f^{i}_{1}}\vf^{i} }$ \ 
of length \ $n_{i}+1$, \ we associate the \
$n_{i}$-dimensional cube \ $C[\Psi^{i}]$ \ -- \ necessarily in \ $\BG$, \
because of the properness \ -- \ and let \ $C^{M}[\Phi]$ \ denote the
subcomplex \  $C[\Psi^{1}]\times C[\Psi^{2}]\times \dotsb\times C[\Psi^{r}]$ \ 
of \ $\BG$ \ -- \ itself a cube. 

If all but \emph{one} of \ $\Psi^{1}$, \ $\Psi^{2},\dotsc$, \ $\Psi^{r}$ \ 
is of length $1$ \ -- \ so that the corresponding factors \ 
$C[\Psi^{i}]$ \ are zero-dimensional \ -- \  we say that \ $C^{M}[\Phi]$ \ is
\emph{indecomposable}. Denote by \ $\Im^{k}$ \ the collection
of all indecomposable $k$-cubes of \ $\BG$ \ (for all maximal chains $\Phi$).

Let \ $\II_{k}$ \ denote the collection of all unordered pairs \ 
$\{C_{i}, C_{j}\}$ \ of indecomposable \ $(k+1)$-cubes of \ $\BG$ \ such
that \ $C_{i,j}\DEF C_{i}\cap C_{j}$ \ is $k$-dimensional.
\end{defn}

\begin{example}\label{epart}{\stepcounter{subsection}}
For \ $\Gamma=L_{4}$ \ as in Example \ref{ecube} above, with  \
$bWL_{4}$ \ as in Figure \ref{fig3}, \ we have three $2$-dimensional
cubes, of which the first two are indecomposable:

\begin{enumerate}
\renewcommand{\labelenumi}{(\alph{enumi})}
\item The left quadrilateral is \ $C_{1}\DEF C[\Psi^{1}]\times C[\Psi^{2}]$, \ 
where \ $\Psi^{1}\DEF\lra{\vi\xra{k}v_{4}}$ \ 
(so \ $C[\Psi^{1}]$ \ is a point) and \ 
$\Psi^{2}=\lra{v_{4} \xra{h} v_{3} \xra{g} v_{2} \xra{f} \vf}$.
\item The upper right square is \ 
$C_{2}\DEF C[\Psi^{3}]\times C[\Psi^{4}]$, \ for \ 
$\Psi^{3}=\lra{\vi \xra{k} v_{4} \xra{h} v_{3} \xra{g} v_{2}}$ \ 
and \ $\Psi^{4}=\lra{v_{2}\xra{f}\vf}$.
\item The lower right quadrilateral decomposes as the product of 
two $1$-dimensional subcubes, corresponding to \ 
$\vi \xra{k} v_{4} \xra{h} v_{3}$ \ and \ $v_{3} \xra{g} v_{2}\xra{f}\vf$, \ 
respectively.
\end{enumerate}

$\II_{2}$ \ consists of the single pair \ $\{C_{1}, C_{2}\}$, \ so
$C_{1,2}=C_{1}\cap C_{2}$ \ is the vertical $1$-cube denoted by \ 
$(f)(g\circ h)(k)$ \ in Figure \ref{fig3}.
\end{example}

\subsection{Conditions for the induction step}
\label{scis}\stepcounter{thm}

Let us assume that \ref{aekr} holds for $k$, so \ $\bA_{k-1}$ \ has been
chosen. Note that because \ $\bA_{k-1}:\sk{k-1}\BG\to\Ta$ \ is a functor,
it extends uniquely to all \emph{decomposable} cubes \ 
$C^{M}[\Phi]=C[\Psi^{1}]\times C[\Psi^{2}]\times \dotsb\times C[\Psi^{r}]$ \
as above, as long as each \ $C[\Psi^{i}]$ \ has dimension 
$\leq k-1$ \ (and it \emph{can} be extended if \ 
$\dim(C[\Psi^{i}])\leq k$ \ for all $i$).

We want to choose a \ $k$-realization \ $\bA_{k}$ \ in such a
way that it has \emph{some} extension to the \ $(k+1)$-skeleton of \ 
$\BG$. \ By the above, it is enough to extend it to the set \ 
$\Im^{k+1}=\{C_{1},\dotsc,C_{\ell}\}$ \ of all indecomposable \ 
$(k+1)$-cubes for $\Gamma$, \ where \ $C_{i}=\WG_{i}$ \ 
for some sublattice \ $\Gamma_{i}$ \ of $\Gamma$ (isomorphic to \ 
$L_{k+2}$ \ of \S \ref{elat}), \ so that \ $\BG_{i}$ \ is $k$-dimensional.

Consider the initial data \ $\A_{i}\DEF\A\rest{\Gamma_{i}}$ \ for each \  
$1\leq i\leq\ell$. \
Clearly,  Assumption \ref{aekr}(b) (for $k$) is needed in
order for each of the higher homotopy operations \
$\llrr{\A_{1}},\dotsc,\llrr{\A_{\ell}}$ \ to even be defined,
(since we need an extension to \ $\BG_{i}$ \ for each $i$); their 
vanishing is \emph{necessary} in order for \ref{aekr}(a) to hold for \
$(k+1)$, \ since it must be possible to extend the new \ $\bA_{k}$ \ 
over the \ $(k+1)$-sekeleton.

However, a \emph{sufficient} condition for \ref{aekr} to hold for \
$(k+1)$ \ is that the operations \ $\llrr{\A_{i}}$ \ vanish
\emph{coherently} \ -- \ which means, essentially, that they can be
made to vanish by choices of realizations which agree on the common
$k$-dimensional faces of the different cubes \ $C_{i}$ \ (and which
equal the given \ $\bA_{k-1}$ \ on the \ $(k-1)$-skeleta of these
faces). \ In order to formulate this condition precisely \ -- \ 
without specifying a global \ $\bA_{k}$ \ for all of \ $\BG$ \ -- \ 
we need the following 

\begin{defn}\label{ddo}\stepcounter{subsection}  
For \ $k\geq 1$, \ let $C$ be a $k$-cube (we have in mind any of the
common facets of two indecomposable \ $(k+1)$-cubes), with boundary \ 
$\partial C=\sk{k-1}C$, \ and let \ $f:X\rtimes\partial C\to Y$ \ be a map.
Note that the space obtained by identifying the boundaries \ 
$\partial C$ \ in two different copies of $C$ is canonically homeomorphic to 
the unreduced suspension \ $\Sigma\partial C$, \ so each pair of extensions \ 
$F_{1},F_{2}:X\rtimes C\to Y$ \ of $f$ defines a map \ 
$D_{F_{1},F_{2}}:X\rtimes\Sigma\partial C\to Y$, \ whose
homotopy class (relative to \ $X\rtimes\partial C$) \ depends 
only on the homotopy classes of \ $F_{1}$ and \ $F_{2}$ \ relative to \ 
$X\rtimes\partial C$. \ This class is called the 
\emph{difference obstruction} for  \ $F_{1}$ and \ $F_{2}$ \ 
\emph{relative to} $f$, and is denoted by
$$
\delta_{f}(F_{1},F_{2})\in
[X\rtimes\Sigma\partial C,Y~~(\text{rel~}X\rtimes\partial C)].
$$ 
\end{defn}

\begin{defn}\label{dobs}\stepcounter{subsection}  
Assume that for each \ $C_{i}=\WG_{i}$ \ ($i=1,\dotsc,\ell$) \ in \
$\Im^{k+1}$ \ we have some extension \ $F_{i}:\WG_{i}\to\Ta$ \ of \ 
$\bA_{k-1}\rest{\sk{k-1}C_{i}}$ (for a $(k-1)$-realization \ 
$\bA_{k-1}$ \ as in \S \ref{aekr}). Note that such an \ $F_{i}$ \
exists if and only if the corresponding higher operation \ 
$\llrr{\A_{i}}$ \ of \S \ref{scis} vanishes.

Then the \emph{obstruction sequence}  \emph{determined by} \ 
$(F_{i})_{i=1}^{\ell}$ \ is the finite sequence 
$$
\Delta(F_{i})_{i=1}^{\ell}\DEF \prod_{\{C_{i}, C_{j}\}\in\II_{k}}\ 
\delta_{\bA_{k-1}}(F_{i}\rest{C_{i,j}},F_{j}\rest{C_{i,j}}).
$$

The obstruction sequence \ $\Delta(F_{i})_{i=1}^{\ell}$ \ is said to 
\emph{vanish} if each of its components \
$\delta_{\bA_{k-1}}(F_{i}\rest{C_{i,j}},F_{j}\rest{C_{i,j}})$ \ is null.
\end{defn}

Note that these classes depend only on the homotopy classes of \ 
$F_{i}\rest{C_{i,j}}$ \ and \ $F_{j}\rest{C_{i,j}}$ \ 
relative to \ $\sk{k-1}C_{i,j}$. \ Evidently:

%
%
\begin{thm}\label{ttwo}\stepcounter{subsection}
Assume given initial data $\A$; if \ $\bA_{k-1}:\sk{k-1}\BG\to\Ta$ \ 
is a \ $(k-1)$-realization which can be extended to a $k$-realization, \ 
then \ $\bA_{k-1}$ \ extends further to a \ $(k+1)$-realization 
if and only if it has some vanishing obstruction sequence.  
\end{thm}

\begin{remark}\label{rcv}\stepcounter{subsection}
This theorem is of little use, if viewed as an obstruction theory for 
successively constructing $k$-realizations of $\A$. It should be
thought of rather as an attempt to make some sense of the claim
that ``an $n$-th order homotopy operation is defined if and only if all
lower order operations vanish coherently.'' The classes in \
$\delta_{\bA_{k-1}}$ \ should be thought of as the obstructions to
\emph{coherence}; the vanishing of the individual operations 
is implicit in Definition \ref{dobs}.
\end{remark}

\begin{example}\label{edo}\stepcounter{subsection}  
In order for the ``long Toda bracket'' \ $\llrr{\A}=\lra{f,g,h,k}$ \ 
of Example \ref{eltb} to be defined, the two (ordinary) Toda brackets \ 
$\lra{f,g,h}$ \ and \ $\lra{g,h,k}$ \ -- \ corresponding to \ $C_{1}$ \ 
and \ $C_{2}$ \ of Example \ref{epart}, respectively \ -- \ must vanish.  
(We may disregard the lower quadrilateral (or triangle), which is
decomposable, and thus represents a ``product'' of ordinary 
nullhomotopies \ $f\circ g\sim\ast$ \ and \ $h\circ k\sim\ast$.)  

By Theorem \ref{ttwo}, the only difference obstruction corresponds 
to \ $C_{1,2}$ \ of Example \ref{epart} \ -- \ the edge \ $(f)(g\circ h)(k)$ \ 
in Figure \ref{fig4}. \ The vanishing of \ $\delta_{\bA_{1}}$ \ 
is clearly equivalent to the choosing of \emph{homotopic} 
nullhomotopies for \ $g\circ h\sim\ast$ \ in defining the two
vanishing elements in \ $\lra{f,g,h}$ \ and \ $\lra{g,h,k}$, \
respectively (relative to fixed choices of $g$ and $h$).
\end{example}

%
%
\sect{Families of polytopes}
\label{cfp}

We have defined higher homotopy operations for any lattice $\Gamma$, as the 
(final) obstruction to rectifying diagrams \ $\A:\Gamma\to\hoT$. \ In this
generality there is very little structure to them. 
However, in most cases of interest $\Gamma$ will have certain ``symmetries''
which will simplify the description of the faces of \ $W\Gamma(\vi,\vf)$ \ 
and ensure that \ $\PP=W\Gamma(\vi,\vf)$ \ will be combinatorially
equivalent to an $n$-dimensional polytope (i.e., convex polyhedron 
in \ $\R^{n}$), \ with basis \ $\BG$ \ equivalent to its boundary \ 
$\partial\PP$, \ which is thus homeomorphic to a sphere. \ 

    This is important because of Remark \ref{rsmash} above, since in 
that situation the corresponding higher operations actually 
take value in appropriate track \emph{groups} (which are abelian, for \ 
$n\geq 3$), \ and the indeterminacy may often be described in terms of 
appropriate cosets, as for the classical Toda bracket 
(cf.\ \cite[Ch.\ I]{TodC}).

\begin{defn}\label{dfop}\stepcounter{subsection}  
A \emph{family of polytopes}  is a sequence \ $\F=(P_{n})_{n=0}^{\infty}$ \ of 
polytopes, starting with \ $P_{0}=\{\pt\}$, \ such that \ $\dim(P_{n})=n$, \ 
and each facet of \ $P_{n}$ \ is isomorphic to some product of lower 
dimensional polytopes from $\F$.
\end{defn}

\begin{examples}\label{efp}\stepcounter{subsection}  
Many familiar examples of polytopes fit into such families:
\begin{enumerate}
\renewcommand{\labelenumi}{\arabic{enumi}.}
\item The family $\Del$ of standard $n$-simplices \ 
$\{\De{n}\}_{n=0}^{\infty}$.
\item The family \ $\Cu$ \ of $n$-cubes \ $\{[0,1]^{n}\}_{n=0}^{\infty}$.
\item The family \ $\Ass=(K_{n})_{n=0}^{\infty}$ \ of \emph{associahedra}, \ 
due to Stasheff (cf.\ \cite{StaH1}).
\item The family \ $\Perm=(\Pe{n})_{n=0}^{\infty}$ \ of \emph{permutohedra}, \ 
apparently due to Schoute at the beginning of the twentieth century 
(cf.\ \cite{SchouA}), where \ $\Pe{n}=\Pe{n}(a_{0},\dotsc, a_{n})$ \ 
is defined to be the convex hull of the \ $(n+1)!$ \ points \ 
$(\sigma(a_{0}),\sigma(a_{1}),\dotsc,\sigma(a_{n}))\in\R^{n+1}$, \
indexed by permutations \ $\sigma\in\Sigma_{n+1}$, \ where \ 
$a_{0},\dotsc, a_{n}$ \ are any distinct real numbers (though we usually take \ 
$(a_{0},\dotsc, a_{n})=(0,1,2,\dotsc,n)$).
\end{enumerate}
\end{examples}

Note that in all of these examples we can find a more economical description 
for $\PP$ than the cubical or simplicial ones described above for general 
$\Gamma$ \ -- \ for example, instead of triangulating the $2$-dimensional 
associahedron using ten $2$-simplices, or five $2$-cubes, we can think of it as
just a single pentagon. Such ``minimal models'' are clearly useful, but there 
appears to be no canonical way to obtain them, in general. 

\subsection{Permutohedra}
\label{sperm}\stepcounter{thm}

The family of permutohedra \ $\Perm=(\Pe{n})_{n=0}^{\infty}$, \ 
which is in some sense the universal family of polytopes,
is associated to the category \ $\N_{+}$ \ of \S \ref{snot}.
The family \ $\Perm$ \ has appeared in a number
of homotopy theoretic contexts \ -- \ in the work of Milgram 
\cite[\S 4]{MilgI}, Stasheff \cite[\S 11]{StaH}, 
Baues \cite[III, (4.5)]{BauG}, and others.

If we let $\Gamma$ denote the full subcategory of \ $\N^{op}_{+}$ \ with 
objects \ $\{\bm,\bz,\dotsc,\bn\}$, \ then we have maps \ 
$d_{i}^{k}:\bk\to\bkm$ \ for \ $0\leq i\leq k\leq n$, \ with \ 
%
%
\setcounter{equation}{\value{thm}}\stepcounter{subsection}
\begin{equation}\label{eseven}
d^{k-1}_{i}d^{k}_{j}=d^{k-1}_{j-1}d^{k}_{i} \ \text{for all} \ 
           0\leq i<j\leq k\leq n
\end{equation}
\setcounter{thm}{\value{equation}}
\noindent (where \ $d^{0}_{0}\DEF\varepsilon$ \ is the augmentation).

The $n$-dimensional permutohedron \ $\Pe{n}$ \ is then a minimal model for 
$\PP$, \ and the simplicial complex \ $\map(\bn,\bm)$ \ which triangulates \ 
$\WG$ \ (see \S \ref{rsv}) is the first barycentric subdivision of the
obvious triangulation of \ $\Pe{n}$. \
The isomorphism may be described explicitly as follows:

Let \ $\Phi^{\Id}$ \ denote the (maximal) $(n+1)$-chain of $\Gamma$, obtained 
by writing the unique map \ $\varphi:\vi=\bn\to\bm=\vf$ \ of $\Gamma$ in 
standard form as \ $d_{0}d_{1}d_{2}\dotsb d_{n}$. \ To the vertex \ 
$v_{\sigma}$ \ ($\sigma\in\Sigma_{n+1}$) \ of \ $\Pe{n}$ \ we then associate 
the vertex indexed by the \ $(n+1)$-chain \ $\Phi^{\sigma}$ \ of $\Gamma$ 
obtained from \ $\Phi^{\Id}$ \ by carrying out the permutation $\sigma$ on 
the maps \ $(d_{0},d_{1},d_{2},\dotsc, d_{n})$ \ and applying \eqref{eseven} 
for each adjacent transposition.

There is an edge $e$ in \ $\Pe{n}$ \ connecting \ $v_{\sigma}$ \ with \ 
$v_{\tau}$ \ whenever $\tau$ is obtained from $\sigma$ by an adjacent 
transposition \ ($i,i+1)$; \ the barycenter of $e$ is indexed by the $n$-chain 
obtained from \ $\Phi^{\sigma}$ \ (or \ $\Phi^{\tau}$) \ by composing
the two maps at the $(i+1)$-st node.  

\begin{remark}\label{rperm}\stepcounter{subsection}
More geometrically, we can think of  $\Gamma$ \ as the lattice 
the $k$-simplices \ $\De{k}$ \ ($0\leq k\leq n$), \ with morphisms
consisting of all possible inclusions of faces.

Alternatively, consider the lattice  \ $\Gamma'$ \ of all subsets of \ 
$\bn=\{0,1,2,\dotsc,n\}$, \ with initial object
$\emptyset$ and terminal object $\bn$; which is equivalent, of course, to
the lattice of all faces of the simplicial complex \ $\De{n}$ \ (with
morphisms again the inclusions). Viewed in this way, \ $\Pe{n}$ \ will be 
called the \emph{face-polyhedron} of \ $\De{n}$. \ 
Many other interesting families of polyhedra can be constructed similarly.

The lattice $\Gamma$ defined above (as well as its opposite \ 
$\Gamma^{op}$) \ is obtained from \ $\Gamma'$ \ by successively identifying 
all subsets of $\bn$ of cardinality $k$ to a single node $\bk$, without 
changing the morphisms. More precisely, if  \ $k<\ell$, \ $B$ is \emph{some} 
subset of $\bn$ of cardinality $\ell$, and \ $A_{1},\dotsc,A_{r}$ \ are the 
$k$-subsets of $B$, then \ 
$\Hom_{\Gamma}(\bk,\bl):=\bigcup_{i=1}^{r}\ \Hom_{\Gamma'}(A_{i},B)$. \ 

As a result, the corresponding cubical categories \ $\WG'$ \ and \ 
$\WG$ \ have homeomorphic mapping spaces \ -- \ and in particular, 
the same total mapping space \ $\PP=\PP_{\Gamma'}\cong\wG$. \ However, for our 
purposes we do need to consider the (distinct) category structures on \ 
$\WG'$ \ and \ $\WG$. 

Of course, this ``condensed'' construction works only because all
$k$-dimensional faces of the simplex are `the same,' (i.e., the group of 
automorphisms act transitively on the set of $k$-dimensional faces, for
all $k$). For a general polyhedron, we expect different faces to correspond to 
different nodes in the lattice of faces $\Gamma$. However, if we are only
interested in a minimal model for \ $\WG$, \ as in the case of the 
permutohedron, the condensed version is just as good.
\end{remark}

\begin{example}\label{esid}\stepcounter{subsection}  
For \ $n=2$ \ the lattice $\Gamma$ of \S \ref{sperm} is described by 
Figure~\ref{fig16}.

\setcounter{figure}{\value{thm}}\stepcounter{subsection}
\begin{figure}[htb]
\begin{center}
%
%
\begin{picture}(220,255)(0,0)
%
%
\put(-8,130){$\vi${\scriptsize{$=\mathbf{2}$}}}
\put(30,130){$\bullet$}
\put(40,134){\vector(1,0){50}}
\put(58,138){$d_{2}$}
\put(96,138){{\scriptsize{$\mathbf{1}$}}}
\put(95,130){$\bullet$}
\put(105,134){\vector(1,0){50}}
\put(118,138){$d_{1}$}
\put(161,138){{\scriptsize{$\mathbf{0}$}}}
\put(160,130){$\bullet$}
\put(170,134){\vector(1,0){50}}
\put(188,138){$d_{0}$}
\put(222,130){$\bullet$}
\put(230,131){{\scriptsize{$\mathbf{-1}=$}}$\vf$}
%
%
\put(87,162){$G_{12}$}
\put(39,139){\vector(1,1){53}}
\put(54,166){$d_{1}$}
\put(95,190){$\bullet$}
\put(96,198){{\scriptsize{$\mathbf{1}$}}}
\put(105,187){\vector(1,-1){50}}
\put(128,166){$d_{1}$}
%
%
\put(157,162){$H$}
\put(104,194){\vector(1,0){50}}
\put(120,198){$d_{0}$}
\put(160,190){$\bullet$}
\put(166,196){{\scriptsize{$\mathbf{0}$}}}
\put(170,187){\vector(1,-1){50}}
\put(192,167){$d_{0}$}
%
%
\put(90,215){$G_{01}$}
\put(36,141){\vector(1,2){57}}
\put(50,195){$d_{0}$}
\put(93,255){$\bullet$}
\put(94,247){{\scriptsize{$\mathbf{1}$}}}
\put(102,255){\vector(1,-1){58}}
\put(133,230){$d_{0}$}
%
%
\put(87,107){$G_{02}$}
\put(36,128){\vector(1,-1){50}}
\put(49,95){$d_{0}$}
\put(90,70){$\bullet$}
\put(91,78){{\scriptsize{$\mathbf{1}$}}}
\put(100,74){\vector(1,0){50}}
\put(111,78){$d_{1}$}
%
%
\put(157,107){$H$}
\put(104,128){\vector(1,-1){50}}
\put(130,105){$d_{0}$}
\put(157,70){$\bullet$}
\put(158,78){{\scriptsize{$\mathbf{0}$}}}
\put(165,77){\vector(1,1){50}}
\put(181,106){$d_{0}$}
%
%
\put(140,45){$H$}
\put(95,65){\vector(1,-1){57}}
\put(104,37){$d_{0}$}
\put(155,5){$\bullet$}
\put(156,13){{\scriptsize{$\mathbf{0}$}}}
\put(163,10){\vector(1,2){58}}
\put(190,55){$d_{0}$}
\end{picture}
\caption{\label{fig16}Lattice for $2$-permutahedron}
\end{center}
\end{figure}
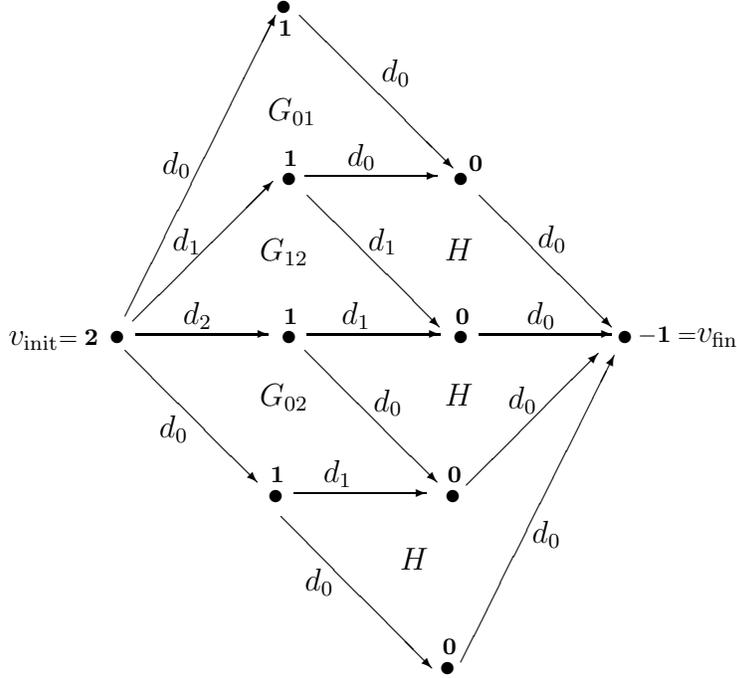
\setcounter{thm}{\value{figure}}
\noindent We have not indicated any of the commuting triangles here,
but only the commuting quadrilaterals which result directly from
applying \eqref{eseven}.  Moreover, the true graph is not planar,
since all edges with the same label and vertices should be identified,
as well as the commuting quadrilaterals labeled $H$. 

The minimal model for \ $\WG$ \ is the $2$-permutohedron \ $\Pe{2}$ in
Figure~\ref{fig6}. 
Here each vertex or edge has been given two different labels \ -- \ 
one corresponding to the description of $\Gamma$ in \S \ref{sperm}, and the
other to the description of \ $\Gamma'$ \ in \S \ref{rperm} (but we
have abbreviated \ $(012,02,\emptyset)$, \ e.g., \ to \ $(012,02)$).

\setcounter{figure}{\value{thm}}\stepcounter{subsection}
\begin{figure}[htb]
\begin{center}
%
%
\begin{picture}(250,190)(20,0)
%
%
\put(155,165){\line(-3,-2){60}}
\put(44,150){\scriptsize $(012,01)=(d_{0}d_{1})d_{2}$}
\put(158,167){\circle*{3}}
\put(109,173){{\scriptsize $(012,01,1)=$}{\small $d_{0}d_{0}d_{2}$}}
%
%
\put(161,165){\line(3,-2){60}}
\put(186,153){\scriptsize $(012,1)=d_{0}(d_{0}d_{2})$}
\put(222,123){\circle*{3}}
\put(226,121){{\scriptsize $(012,12,1)=$}{\small $d_{0}d_{1}d_{0}$}}
%
%
\put(222,120){\line(0,-1){60}}
\put(224,92){\scriptsize $(012,12)=(d_{0}d_{1}) d_{0}$}
\put(222,58){\circle*{3}}
\put(225,56){{\scriptsize $(012,12,2)=$}{\small $d_{0}d_{0}d_{0}$}}
%
%
\put(221,59){\line(-3,-2){60}}
\put(192,32){\scriptsize $(012,2)=d_{0}(d_{0}d_{1})$}
\put(108,6){{\scriptsize $(012,02,2)=\,$}{\small $d_{0}d_{0}d_{1}$}}
%
%
\put(94,58){\circle*{3}}
\put(5,56){{\scriptsize $(012,02,0)=$}{\small $d_{0}d_{1}d_{1}$}}
\put(154,18){\line(-3,2){60}}
\put(158,18){\circle*{3}}
\put(45,30){\scriptsize $(012,02)=(d_{0}d_{1})d_{1}$}
%
%
\put(94,123){\circle*{3}}
\put(5,121){{\scriptsize $(012,01,0)=$}{\small $d_{0}d_{1}d_{2}$}}
\put(94,120){\line(0,-1){60}}
\put(14,90){\scriptsize $(012,0)=d_{0}(d_{1}d_{2})$}
\end{picture}
\caption{\label{fig6}Two descriptions of \protect{$\Pe{2}$}}
\end{center}
\end{figure}
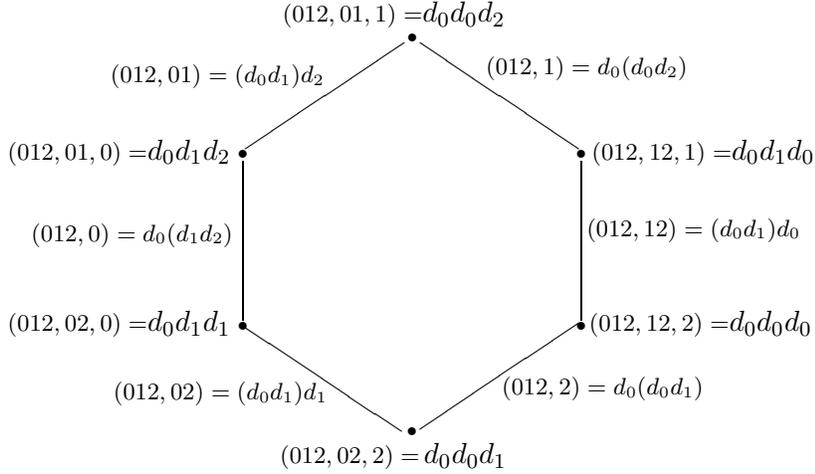
\setcounter{thm}{\value{figure}}
\end{example}

Similarly, for \ $n=3$ \ we obtain the $3$-permutahedron \ $\Pe{3}$ \ 
(with maximal chain \ $d_{0}d_{1}d_{2}d_{3}$) in Figure~\ref{fig8}.
Note that we have indicated only the labels of \S \ref{sperm}, where 
the chain \ $d_{i_{0}}d_{i_{1}}d_{i_{2}}d_{i_{3}}$, \ for example,
is represented by \ $(i_{0},i_{1},i_{2},i_{3})$.

\setcounter{figure}{\value{thm}}\stepcounter{subsection}
\begin{figure}[htb]
\begin{center}
%
%
\begin{picture}(200,190)(30,-10)
%
%
\put(40,72){\circle*{5}}
\put(0,69){\scriptsize $(0,1,2,2)$}
\put(41,69){\line(1,-3){10}}
\put(52,36){\circle*{5}}
\put(13,33){\scriptsize $(0,1,1,2)$}
\put(54,36){\line(1,0){33}}
\put(88,36){\circle*{5}}
\put(57,42){\scriptsize $(0,1,1,1)$}
\put(89,37){\line(2,3){22}}
\put(112,72){\circle*{5}}
\put(74,70){\scriptsize $(0,1,2,1)$}
\put(41,74){\line(1,3){10}}
\put(52,108){\circle*{5}}
\put(13,105){\scriptsize $(0,1,2,3)$}
\put(55,108){\line(1,0){45}}
\put(100,108){\circle*{5}}
\put(62,111){\scriptsize $(0,1,1,3)$}
\put(101,105){\line(1,-3){10}}
%
%
\put(54,34){\line(4,-3){47}}
\put(100,0){\circle*{5}}
\put(58,-4){\scriptsize $(0,0,1,2)$}
\put(103,0){\line(1,0){21}}
\put(124,0){\circle*{5}}
\put(128,-8){\scriptsize $(0,0,1,1)$}
\put(91,32){\line(1,-1){33}}
%
%
\put(128,0){\line(4,1){45}}
\put(172,12){\circle*{5}}
\put(177,6){\scriptsize $(0,0,0,1)$}
\put(173,15){\line(1,4){8}}
\put(184,60){\circle*{5}}
\put(170,50){\scriptsize $(0,1,0,1)$}
\put(115,73){\line(3,1){34}}
\put(148,84){\circle*{5}}
\put(127,71){\scriptsize $(0,0,2,1)$}
\put(151,82){\line(3,-2){32}}
%
%
\put(174,14){\line(1,1){45}}
\put(220,60){\circle*{5}}
\put(225,57){\scriptsize $(0,0,0,0)$}
\put(221,64){\line(1,4){10}}
\put(186,62){\line(1,1){10}}
\put(205,81){\line(1,1){23}}
\put(232,108){\circle*{5}}
\put(235,105){\scriptsize $(0,1,0,0)$}
%
%
\put(136,120){\circle*{5}}
\put(100,124){\scriptsize $(0,0,1,3)$}
\put(137,117){\line(1,-3){6}}
\put(144,96){\line(1,-3){2}}
\put(138,122){\line(3,4){33}}
\put(172,168){\circle*{5}}
\put(177,171){\scriptsize $(0,0,0,3)$}
\put(175,167){\line(2,-1){44}}
\put(220,144){\circle*{5}}
\put(225,141){\scriptsize $(0,0,2,0)$}
\put(221,141){\line(1,-3){10}}
%
%
\put(103,109){\line(3,1){32}}
\put(54,111){\line(2,3){21}}
\put(76,144){\circle*{5}}
\put(37,141){\scriptsize $(0,0,2,3)$}
\put(79,145){\line(2,1){45}}
\put(124,168){\circle*{5}}
\put(112,173){\scriptsize $(0,1,0,3)$}
\put(127,168){\line(1,0){45}}
%
%
\multiput(77,140)(1,-4){12}{\circle*{.5}}
\put(88,96){\circle*{3}}
\put(60,87){\scriptsize $(0,0,2,2)$}
\multiput(91,95)(3,-1){12}{\circle*{.5}}
\put(124,84){\circle*{3}}
\put(112,90){\scriptsize $(0,1,0,2)$}
\multiput(127,86)(3,2){2}{\circle*{.5}}
\multiput(145,98)(3,2){6}{\circle*{.5}}
\put(160,108){\circle*{3}}
\put(165,105){\scriptsize $(0,1,1,0)$}
\multiput(161,111)(1,3){12}{\circle*{.5}}
\put(172,144){\circle*{3}}
\put(171,134){\scriptsize $(0,1,2,0)$}
\multiput(168,146)(-4,2){12}{\circle*{.5}}
%
%
\multiput(44,74)(4,2){6}{\circle*{.5}}
\put(84,94){\circle*{.5}}
%
%
\multiput(123,81)(1,-3){12}{\circle*{.5}}
\put(136,48){\circle*{3}}
\put(135,40){\scriptsize $(0,0,0,2)$}
\multiput(133,46)(-3,-4){12}{\circle*{.5}}
%
%
\multiput(141,50)(3,2){12}{\circle*{.5}}
\put(172,72){\circle*{3}}
\put(174,74){\scriptsize $(0,0,1,0)$}
\multiput(175,71)(4,-1){12}{\circle*{.5}}
%
%
\multiput(162,105)(1,-3){12}{\circle*{.5}}
\multiput(176,144)(4,0){12}{\circle*{.5}}
\end{picture}
\caption{\label{fig8}\protect{A depiction of $\Pe{3}$}}
\end{center}
\end{figure}
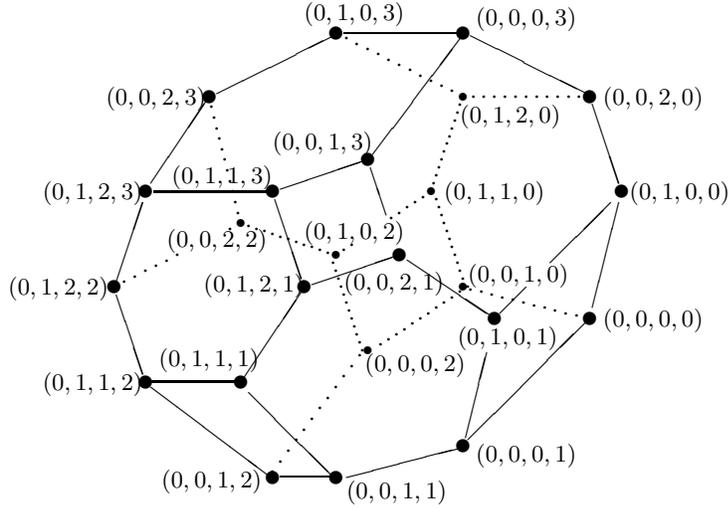
\setcounter{thm}{\value{figure}}

\begin{remark}\label{rhperm}\stepcounter{subsection}
The higher homotopy operations associated to the family \ $\Perm$ \ were
described explicitly in \cite[\S 5]{BlaHH}. In fact, it was this  case
which motivated the search for a general definition of higher
operations, which culminated in this paper.
\end{remark}

\begin{example}\label{ehperm}\stepcounter{subsection}
Let us describe explicitly the higher homotopy operation associated to the
lattice $\Gamma$ of \S \ref{sperm} for $n=2$.

A homotopy $\Gamma$-diagram is clearly given by maps \ 
$f_{0},f_{1},f_{2} : X \to Y$, \ $g_{0},g_{1} : Y \to Z$, \ and \ 
$h : Z \to W$, \ as well as maps \ $A: Y \to W$, \ 
$B_{01}, B_{02}, B_{12} : X \to Z$ and $C: X \to W$, \ 
such that \ $[h g_{0}] = [h g_{1}] = [A]$, \ 
$[g_{0} f_{0}] = [g_{0} f_{1}] = [B_{01}]$, \ 
$[g_{1} f_{1}] = [g_{1} f_{2}] = [B_{12}]$, \ 
$[g_{0} f_{2}]= [g_{1} f_{1}] = [B_{02}]$ \ and $[h g_{0} f_{0}] = [C]$.

One readily verifies that the maps \ $A$, \ $B_{01}$, \ $B_{02}$, \ 
$B_{12}$ \ and $C$ are irrelevant to the homotopy operation, and may
disregarded \ -- \ this would correspond to taking the `minimal model'
of the hexagon \ $\Pe2$ \ in Figure~\ref{fig6}, instead of its
decomposition into six cubes using the cubical structure of the 
$W$-construction.

Complete data are then given by explicit homotopies \ 
$H: h g_{0} \sim h g_{1}$, \ $G_{01}: g_{0} f_{0} \sim g_{0} f_{1}$, \  
$G_{12}: g_{1} f_{1} \sim g_{1} f_{2}$ \ and \ 
$G_{02}: g_{0} f_{2} \sim g_{1} f_{1}$, \ which can be organized into 
a map \ $X\rtimes S^{1}\to W$.

The corresponding higher homotopy operation is then the subset \ 
$\lra{h,g_{0},g_{1},f_{0},f_{1},f_{2}}\subset[X\rtimes S^{1},W]$ \ 
formed by all complete data, as in Figure \ref{fig7} (where we have
also indicated by dotted lines the cubical decomposition provided by
the maps \  $A$, \ $B_{01}$, \ $B_{02}$, \ and \ $B_{12}$).

\setcounter{figure}{\value{thm}}\stepcounter{subsection}
\begin{figure}[htb]
\begin{center}
%
%
\begin{picture}(250,170)(20,0)
\put(155,90){
\makebox(0,0){\circle*{2}}
\multiput(0,0)(4,0){16}{\makebox(0,0){\hskip -.4mm $.$}}
\multiput(0,0)(-4,0){16}{\makebox(0,0){\hskip -.4mm $.$}}
\multiput(0,0)(2,3.33){17}{\makebox(0,0){\hskip -.4mm $.$}}
\multiput(0,0)(-2,3.33){17}{\makebox(0,0){\hskip -.4mm $.$}}
\multiput(0,0)(2,-3.33){16}{\makebox(0,0){\hskip -.4mm $.$}}
\multiput(0,0)(-2,-3.33){16}{\makebox(0,0){\hskip -.4mm $.$}}
}
%
%
\put(155,165){\line(-3,-2){60}}
\put(100,150){$H f_2$}
\put(158,167){\circle*{3}}
\put(145,173){{$h g_{0} f_{2}$}}
%
%
\put(161,165){\line(3,-2){60}}
\put(190,153){$h G_{02}$}
\put(222,123){\circle*{3}}
\put(230,121){{$h g_1 f_ 0$}}
%
%
\put(222,120){\line(0,-1){60}}
\put(230,92){$H f_0$}
\put(222,58){\circle*{3}}
\put(230,56){{$h g_0 f_0$}}
%
%
\put(221,59){\line(-3,-2){60}}
\put(200,32){$h G_{01}$}
\put(150,6){{$h g_0 f_1$}}
%
%
\put(94,58){\circle*{3}}
\put(60,56){{$h g_1 f_1$}}
\put(154,18){\line(-3,2){60}}
\put(158,18){\circle*{3}}
\put(95,30){$H f_1$}
%
%
\put(94,123){\circle*{3}}
\put(60,121){{$h g_1 f_2$}}
\put(94,120){\line(0,-1){60}}
\put(60,90){$h G_{12}$}
\end{picture}
\caption{\label{fig7}Complete data for the operation 
\protect{$\lra{h,g_{0},g_{1},f_{0},f_{1},f_{2}}$}}
\end{center}
\end{figure}
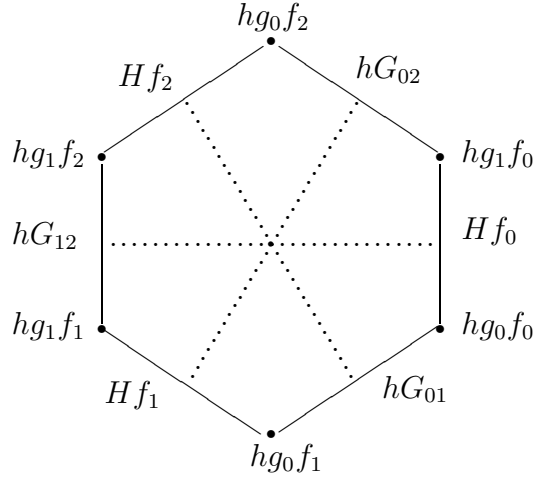
\setcounter{thm}{\value{figure}}
\noindent To describe higher homotopy operations related to other 
polyhedra of this section is equally straightforward and we may safely 
leave these calculations to the reader as an exercise.
\end{example}

\subsection{Simplices}
\label{ssimp}\stepcounter{thm}

Many other families of polytopes are obtained from \ $\Perm$ \ by considering 
\emph{relative} higher homotopy operations (\S \ref{sro}), which has
the effect of collapsing some of the cubes (or simplices) of \ $\WG$. \ 

Thus, the family \ $\Del=\{\De{n}\}_{n=0}^{\infty}$ \ of simplices arises 
from the same diagram \ $\N_{+}$ \ as \ $\Perm$, \ in the case where all the 
relations of \eqref{eseven} are assumed to hold on the nose, \emph{except} for 
the single relation \ $\varepsilon d_{0}=\varepsilon d_{1}$. \ In this case the
permutohedron \ $\Pe{n}$ \ (corresponding to the subcategory $\Gamma$ of \
$\N_{+}^{\op}$ \ with objects \ $\{\bm,\bz,\dotsc,\bn\}$) \ collapses
canonically to an $n$-simplex. 

For example, Figure \ref{fig10} shows the $3$-simplex corresponding to the
possible decompositions of \ $\varepsilon d_{0}d_{1}d_{2}$.

\setcounter{figure}{\value{thm}}\stepcounter{subsection}
\begin{figure}[htbp]
\begin{center}
%
%
\begin{picture}(250,200)(0,-10)
%
%
\put(0,170){\scriptsize{right face} $0$}
\put(200,170){\scriptsize{left face} $1$}
\put(0,5){\scriptsize{back face} $2$}
\put(200,5){\scriptsize{bottom face} $3$}
%
%
\put(90,175){\circle*{3}}
\put(87,172){\vector(-1,-1){59}}
\put(93,173){\scriptsize{$002$}}
\put(115,170){$=$}
\put(95,165){\circle*{3}}
\put(98,163){\scriptsize{$010$}}
\put(92,162){\vector(-1,-1){58}}
\put(62,141){$=$}
\put(45,145){$02$}
\put(75,135){$10$}
\put(110,150){\circle*{3}}
\put(110,145){\vector(0,-1){113}}
\put(103,153){\scriptsize{$000$}}
\put(130,150){\circle*{3}}
\put(130,145){\vector(0,-1){113}}
\put(124,153){\scriptsize{$001$}}
\put(115,112){$=$}
\put(96,112){$00$}
\put(134,112){$01$}
\put(145,165){\circle*{3}}
\put(148,162){\vector(1,-1){58}}
\put(130,163){\scriptsize{$011$}}
\put(155,175){\circle*{3}}
\put(140,173){\scriptsize{$012$}}
\put(158,172){\vector(1,-1){59}}
\put(173,140){$=$}
\put(186,145){$12$}
\put(155,133){$11$}

%
%
\put(25,110){\circle*{3}}
\put(10,107){\scriptsize{$102$}}
\put(30,100){\circle*{3}}
\put(33,97){\scriptsize{$110$}}
\put(16,90){\circle*{3}}
\put(0,88){\scriptsize{$022$}}
\multiput(20,90)(5,0){18}{\line(1,0){3}}
\multiput(113,90)(5,0){3}{\line(1,0){3}}
\multiput(133,90)(5,0){18}{\line(1,0){3}}
\put(224,90){\vector(1,0){2}}
\put(16,80){\circle*{3}}
\put(0,78){\scriptsize{$023$}}
\multiput(20,80)(5,0){18}{\line(1,0){3}}
\multiput(113,80)(5,0){3}{\line(1,0){3}}
\multiput(133,80)(5,0){18}{\line(1,0){3}}
\put(224,80){\vector(1,0){2}}
\put(70,82){$=$}
\put(70,93){$22$}
\put(70,69){$23$}
\put(27,67){\circle*{3}}
\put(11,65){\scriptsize{$103$}}
\put(39,72){\circle*{3}}
\put(42,69){\scriptsize{$120$}}
%
%
\put(220,110){\circle*{3}}
\put(223,107){\scriptsize{$112$}}
\put(206,100){\circle*{3}}
\put(191,97){\scriptsize{$111$}}
\put(230,90){\circle*{3}}
\put(233,88){\scriptsize{$122$}}
\put(230,80){\circle*{3}}
\put(233,78){\scriptsize{$123$}}
\put(219,67){\circle*{3}}
\put(222,65){\scriptsize{$113$}}
\put(204,74){\circle*{3}}
\put(189,71){\scriptsize{$121$}}
%
%
\put(90,5){\circle*{3}}
\put(87,7){\vector(-1,1){57}}
\put(93,3){\scriptsize{$003$}}
\put(115,7){$=$}
\put(95,15){\circle*{3}}
\put(92,17){\vector(-1,1){49}}
\put(98,13){\scriptsize{$020$}}
\put(60,35){$=$}
\put(42,30){$03$}
\put(73,39){$20$}
\put(110,30){\circle*{3}}
\put(103,22){\scriptsize{$100$}}
\put(130,30){\circle*{3}}
\put(124,22){\scriptsize{$101$}}
\put(145,15){\circle*{3}}
\put(148,17){\vector(1,1){50}}
\put(130,13){\scriptsize{$021$}}
\put(155,5){\circle*{3}}
\put(158,7){\vector(1,1){57}}
\put(140,3){\scriptsize{$013$}}
\put(175,35){$=$}
\put(190,28){$13$}
\put(155,39){$22$}
\end{picture}
\end{center}
\caption{\label{fig10}The $3$-simplex}
\end{figure}
\setcounter{thm}{\value{figure}}

\begin{remark}\label{rhsimp}\stepcounter{subsection}
The higher homotopy operations associated to the family \ $\Del$ \  were
described explicitly in \cite[\S 4]{BlaHO}.
\end{remark}

\subsection{Associahedra}
\label{sass}\stepcounter{thm}

Similarly, the family \ $\Ass=(K_{n})_{n=0}^{\infty}$ \ of associahedra arises 
from a less drastic relative version of the diagram \ $\N_{+}$, \ using
a procedure defined by Andy Tonks for collapsing certain faces of the 
permutohedron; this can be defined explicitly in terms of the lattice \ 
$\N_{+}$. \ We shall just give two examples, refering the reader to 
\cite{TonkR} for further details.

\begin{example}\label{eda}\stepcounter{subsection}  
If we think of the $2$-permutohedron as the face polyhedron for \
$\De{2}$ \ (\S \ref{rperm}), then the edge \ $(012,02)$  \ in
Figure \ref{fig6} should be collapsed to yield the $2$-associahedron
in Figure~\ref{fig11}, where the dotted edge is collapsed.

\setcounter{figure}{\value{thm}}\stepcounter{subsection}
\begin{figure}[htb]
\begin{center}
%
%
%
\begin{picture}(120,150)(55,-10)
\unitlength.5mm
%
%
\put(95,85){\line(-3,-2){30}}
\put(47,79){\scriptsize $(012,01)$}
\put(97,87){\circle*{3}}
\put(72,93){\small $(012,01,1)$}
%
%
\put(100,85){\line(3,-2){30}}
\put(116,78){\scriptsize $(012,1)$}
\put(132,63){\circle*{3}}
\put(135,61){\small $(012,12,1)$}
%
%
\put(132,60){\line(0,-1){30}}
\put(134,42){\scriptsize $(012,12)$}
\put(132,28){\circle*{3}}
\put(135,24){\small $(012,12,2)$}
%
%
\put(130,26){\line(-3,-1){40}}
\put(110,11){\scriptsize $(012,2)$}
\put(88,11){\circle*{3}}
\put(68,0){\small $(012,02,2)$}
%
%
\put(74,22){\circle*{3}}
\put(18,22){\small $(012,02,0)$}
\multiput(76,19)(3,-2){4}{\circle*{.5}}
\put(43,11){\scriptsize $(012,02)$}
%
%
\put(64,63){\circle*{3}}
\put(9,60){\small $(012,01,0)$}
\put(64,61){\line(1,-4){9}}
\put(34,41){\scriptsize $(012,0)$}
\end{picture}
\caption{\label{fig11}The $2$-associahedron}
\end{center}
\end{figure}
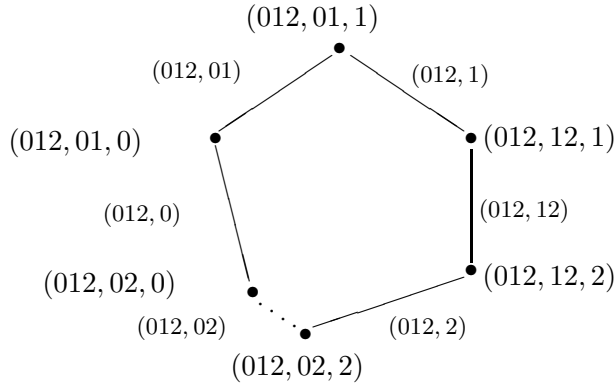
\setcounter{thm}{\value{figure}}

Similarly, in the $3$-permutahedron of Figure \ref{fig8} with the
labeling of \S \ref{rperm}, we collapse each of the edges \ 

\begin{center}
$(0123,012,02)$,  
$(0123,023,02)$,  
$(0123,023,03)$,  
$(0123,013,03)$,  
$(0123,023,2)$,
\\  
$(0123,023,3)$,  
$(0123,123,13)$,  
$(0123,013,13)$,  
$(0123,013,1)$ and  
$(0123,013,0)$,  
\end{center}

\noindent of \ $\Pe{3}$ \ to a point, so that each of the $2$-dimensional 
faces \ 

\begin{center}
$(0123,13)$, \ 
$(0123,013)$, \ 
$(0123,02)$,  \ 
$(0123,03)$ and
$(0123,023)$ \ 
\end{center}

\noindent of $\Pe{3}$ \ is collapsed to an edge in the resulting 
$3$-associahedron.

One can also describe the $n$-associahedron directly in terms of 
bracketing on \ $(n+2)$ \ symbols (cf.\ \cite[\S 2]{StaH1}) \ -- \ 
equivalently, in terms of labelled trees (cf.\ \cite[\S 11]{StaH}) \ -- \ 
or as a truncated $n$-simplex (cf.\ \cite[\S 5.1]{SSterQ});  but these 
descriptions do not fit into our framework of lattices.  
\end{example}

\begin{remark}\label{rhass}\stepcounter{subsection}
The higher homotopy operations associated to the family \ $\Ass$ \
have not been described in this language, but implicitly they
motivated Stasheff's original definition in \cite{StaH1}. 
\end{remark}

\begin{example}
Let $\Gamma$ be as in~\ref{sperm}. Consider the relative case when all
equations of~(\ref{eseven}) are strictly preserved except for
$$
\hskip 1cm
d_0^{k-1} d_1^k = d_0^{k-1}d_0^k,\ 1 \leq k \leq n.
$$
\noindent The corresponding family of polyhedra is the sequence \
$\Cu$ \ of cubes presented as quotients of the permutohedra, as 
illustrated for \ $n=2$ \ in Figure~\ref{fig12}.

\setcounter{figure}{\value{thm}}\stepcounter{subsection}
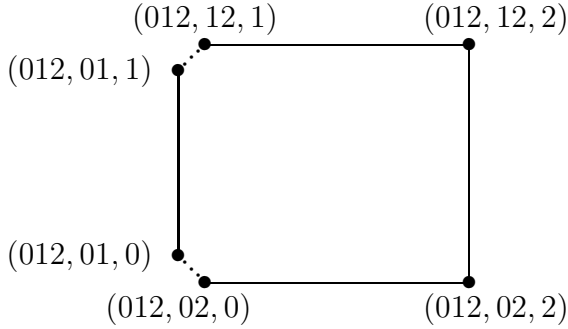
\begin{figure}[htb]
\begin{center}
%
%
{\unitlength=1pt
\thinlines
\begin{picture}(130,120)(0,0)
\put(20,110){\makebox(0,0){$(012,12,1)$}}
\put(0,90){\makebox(0,0)[r]{$(012,01,1)$}}
\put(0,20){\makebox(0,0)[r]{$(012,01,0)$}}
\put(10,0){\makebox(0,0){$(012,02,0)$}}
\put(130,0){\makebox(0,0){$(012,02,2)$}}
\put(130,110){\makebox(0,0){$(012,12,2)$}}
\put(120,100){\makebox(0,0){$\bullet$}}
\put(120,10){\makebox(0,0){$\bullet$}}
\put(20,10){\makebox(0,0){$\bullet$}}
\multiput(10,20)(2.5,-2.5){5}{\makebox(0,0){$\cdot$}}
\put(10,20){\makebox(0,0){$\bullet$}}
\multiput(10,90)(2.5,2.5){5}{\makebox(0,0){$\cdot$}}
\put(10,90){\makebox(0,0){$\bullet$}}
\put(20,100){\makebox(0,0){$\bullet$}}
\put(10,20){\line(0,1){70}}
\put(120,10){\line(-1,0){100}}
\put(120,100){\line(0,-1){90}}
\put(20,100){\line(1,0){100}}
\end{picture}}

\end{center}
\caption{\label{fig12}The $2$-cube as a quotient of the $2$-permutahedron}
\end{figure}
\setcounter{thm}{\value{figure}}
\end{example}

\begin{example}
Consider two copies \ $\Gamma'$ \ and \ $\Gamma''$ \ of the lattice
from~\ref{sperm}, with objects \ 
$\{\mathbf{-1}',\mathbf{0}',\dotsc,\mathbf{n}'\}$ \ and \ 
$\{\mathbf{-1}'',\mathbf{0}'',\dotsc,\mathbf{n}''\}$ \ respectively.
Let $\Gamma$ be the lattice obtained from the disjoint union \ 
$\Gamma' \sqcup \Gamma''$ \ by adding arrows \ 
$f_i:\mathbf{i}' \to \mathbf{i}''$ \ for \ $-1 \leq i \leq n$ \ 
satisfying
%
%
\setcounter{equation}{\value{thm}}\stepcounter{subsection}
\begin{equation}\label{esix}
f_{k-1}d^k_j=d^k_j f_k,\text{\ for\ }0\leq k\leq n \text{\ and\ }0\leq j\leq k
\end{equation}
\setcounter{thm}{\value{equation}}
\noindent See Figure~\ref{fig13}.

\setcounter{figure}{\value{thm}}\stepcounter{subsection}
\begin{figure}[htb]
%
%
\begin{center}
{
\unitlength=1.3pt
\begin{picture}(90,170)(0,0)
\put(40,166){\makebox(0,0){$f_{n}$}}
\put(20,158.50){\vector(1,0){50}}
\put(80,160){\makebox(0,0){$\mathbf{n}''$}}
\put(10,160){\makebox(0,0){$\mathbf{n}'$}}

\put(83,70){\vector(0,-1){20}}
\put(77.50,70){\vector(0,-1){20}}
\put(13,70){\vector(0,-1){20}}
\put(7,70){\vector(0,-1){20}}

\put(85,110){\vector(0,-1){20}}
\put(80,110){\vector(0,-1){20}}
\put(75,110){\vector(0,-1){20}}
\put(15,110){\vector(0,-1){20}}
\put(5,110){\vector(0,-1){20}}
\put(10,110){\vector(0,-1){20}}

\put(80,30){\vector(0,-1){20}}
\put(10,30){\vector(0,-1){20}}
\put(40,130){\makebox(0,0){$f_{2}$}}
\put(40,90){\makebox(0,0){$f_{1}$}}
\put(40,50){\makebox(0,0){$f_{0}$}}
\put(40,10){\makebox(0,0){$f_{-1}$}}
\put(20,120){\vector(1,0){50}}
\put(20,80){\vector(1,0){50}}
\put(20,40){\vector(1,0){50}}
\put(20,0){\vector(1,0){50}}
\put(80,140){\makebox(0,0){$\vdots$}}
\put(10,140){\makebox(0,0){$\vdots$}}
\put(80,120){\makebox(0,0){$\mathbf{2}''$}}
\put(80,80){\makebox(0,0){$\mathbf{1}''$}}
\put(80,40){\makebox(0,0){$\mathbf{0}''$}}
\put(80,0){\makebox(0,0){$\mathbf{-1}''$}}
\put(10,120){\makebox(0,0){$\mathbf{2}'$}}
\put(10,80){\makebox(0,0){1$\mathbf{1}'$}}
\put(10,40){\makebox(0,0){$\mathbf{0}'$}}
\put(10,0){\makebox(0,0){$\mathbf{-1}'$}}
\end{picture}}
\end{center}
\caption{\label{fig13}The category $\Gamma$}

\end{figure}
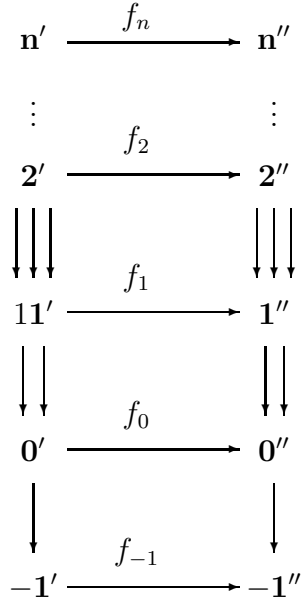
\setcounter{thm}{\value{figure}}

Obviously, a $\Gamma$-diagram consists of two truncated
$\Delta$-simplicial spaces and their $\Delta$-simplicial homomorphism.
  
Consider the relative case when all simplicial identities are strictly
satisfied and also all the identities of \eqref{esix} are strict, except
for the case when $j=k$. The relevant polyhedron is now the
$(n+1)$-simplex obtained as a quotient of the permutahedron with
$(n+2)!$-vertices indexed by
\begin{center}
$d^0_{i_0}d^1_{i_1}\dotsb d^n_{i_n}f_n$, 
$d^0_{i_0}d^1_{i_1}\dotsb d^{n-1}_{i_{n-1}} f_{n-1} d^n_{i_n}$,\ldots
$d^0_{i_0}f_0d^1_{i_1}\dotsb d^n_{i_n}$ and
$f_{-1}d^0_{i_0}d^1_{i_1}\dotsb d^n_{i_n}$,
\end{center}
\noindent for \ $0 \leq i_j \leq j$, \ $0 \leq j \leq n$. \ This is 
illustrated, for \ $n=1$, \ in Figure~\ref{fig14}.

\setcounter{figure}{\value{thm}}\stepcounter{subsection}
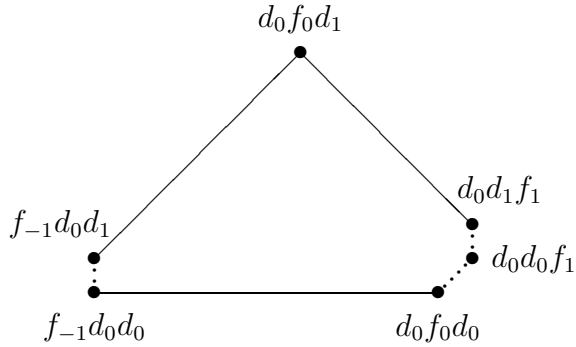
\begin{figure}[htbp]
%
%
\begin{center}
{ \unitlength=1.3pt
\begin{picture}(150,90)(0,0)
\put(110,0){\makebox(0,0){$d_0f_0d_0$}}
\put(138,20){\makebox(0,0){$d_0d_0f_1$}}
\put(128,40){\makebox(0,0){$d_0d_1f_1$}}
\put(70,90){\makebox(0,0){$d_0f_0d_1$}}
\put(0,30){\makebox(0,0){$f_{-1}d_0d_1$}}
\put(10,0){\makebox(0,0){$f_{-1}d_0d_0$}}
\multiput(120,20)(0,3.2){4}{\makebox(0,0){$\cdot$}}
\multiput(110,10)(2,2){5}{\makebox(0,0){$\cdot$}}
\multiput(10,10)(0,3.2){4}{\makebox(0,0){$\cdot$}}
\put(120,20){\makebox(0,0){$\bullet$}}
\put(120,30){\makebox(0,0){$\bullet$}}
\put(110,10){\makebox(0,0){$\bullet$}}
\put(10,10){\makebox(0,0){$\bullet$}}
\put(10,20){\makebox(0,0){$\bullet$}}
\put(70,80){\makebox(0,0){$\bullet$}}
\put(10,10){\line(1,0){100}}
\put(70,80){\line(1,-1){50}}
\put(10,20){\line(1,1){60}}
\end{picture} }
\end{center}
\caption{\label{fig14}The $2$-simplex as a quotient of the $2$-permutahedron}
\end{figure}
\setcounter{thm}{\value{figure}}
\end{example}

%
%
\sect{Massey products}
\label{chmp}

We now show how Massey products fit into our general framework. 

\subsection{Higher Whitehead products}
\label{shwp}\stepcounter{thm}

Because it fits better with our original definition, we start with the Lie 
analogue, sometimes called higher-order Whitehead products:

Any three elements in \ $\pis X$ \ determine a map of the form \ 
$F:S^{r}\vee S^{s}\vee S^{t}\to X$; \ if two of their pairwise  Whitehead 
products vanish, $F$ fits into the diagram in Figure~\ref{fig15},
where \ $w_{r,s}:S^{r+s-1}\to S^{r}\vee S^{s}$ \ is the Whitehead 
product map (and the unmarked maps are the obvious inclusions or projections).

\setcounter{figure}{\value{thm}}\stepcounter{subsection}
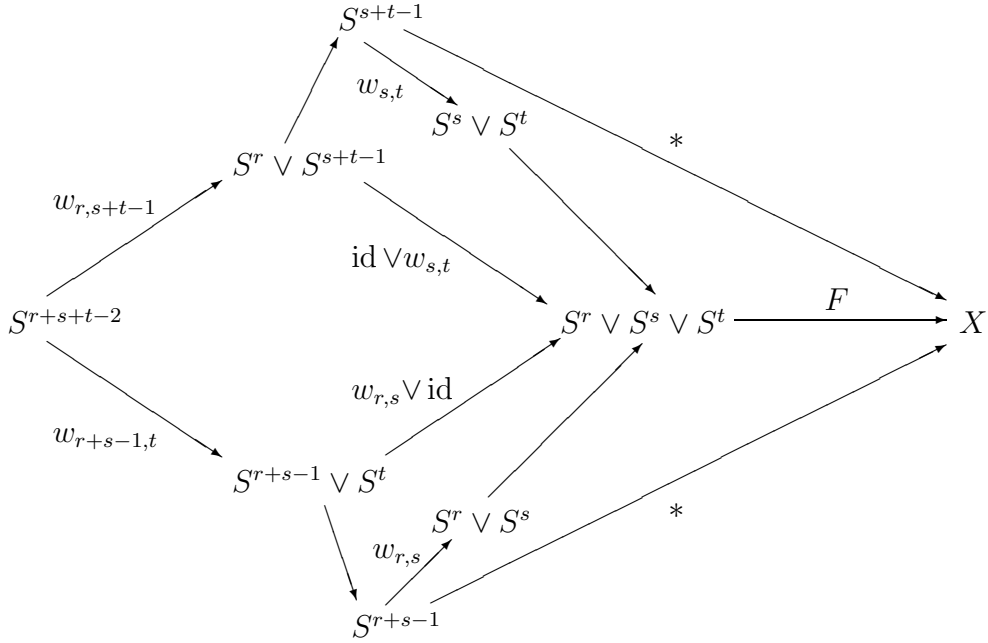
\begin{figure}[htbp]
\begin{center}
%
%
\begin{picture}(367,250)(0,40)
%
%
\put(0,160){$S^{r+s+t-2}$}
\put(15,174){\vector(3,2){66}}
\put(17,207){$w_{r,s+t-1}$}
\put(85,220){$S^{r}\vee S^{s+t-1}$}
\put(15,157){\vector(3,-2){66}}
\put(17,119){$w_{r+s-1,t}$}
\put(85,100){$S^{r+s-1}\vee S^{t}$}
%
%
\put(135,217){\vector(3,-2){69}}
\put(130,185){id\,$\vee w_{s,t}$}
\put(209,160){$S^{r}\vee S^{s}\vee S^{t}$}
\put(143,114){\vector(3,2){66}}
\put(130,135){$w_{r,s}\vee$\,id}
%
%
\put(105,232){\vector(1,2){20}}
\put(125,275){$S^{s+t-1}$}
\put(135,270){\vector(3,-2){35}}
\put(132,251){$w_{s,t}$}
\put(160,235){$S^{s}\vee S^{t}$}
\put(190,230){\vector(1,-1){55}}
%
\put(150,275){\vector(2,-1){205}}
\put(250,230){$\ast$}
%
%
\put(120,95){\vector(1,-3){12}}
\put(130,45){$S^{r+s-1}$}
\put(143,57){\vector(1,1){25}}
\put(138,75){$w_{r,s}$}
\put(160,85){$S^{r}\vee S^{s}$}
\put(182,98){\vector(1,1){58}}
%
\put(160,58){\vector(2,1){195}}
\put(250,90){$\ast$}
%
%
\put(275,165){\vector(1,0){80}}
\put(309,169){$F$}
\put(360,160){$X$}
\end{picture}
\caption{\label{fig15}The diagram defining triple Whitehead products}
\end{center}
\end{figure}
\setcounter{thm}{\value{figure}}

If we think of this as a functor \ $\A:\Gamma\to ho\Ta$ \ for the 
corresponding lattice $\Gamma$, and try to rectify it (relative to the
specified null maps \ -- \ or even relative to the subdiagram involving only 
maps between wedges of spheres), the first (and only) obstruction \ 
$\llrr{\A}\subseteq \pi_{r+s+t-1}X$ \ is called a 
\emph{secondary Whitehead product}. One has $n$-th order
analgues defined for sets of \ $n+1$ \ elements in \ $\pis X$ \ for which all
lower-order Whitehead products vanish coherently (see \cite{GPorW}).

Of course, all this is valid not only in the category of topological spaces,
but also in the category of differential graded Lie algebras (DGLs); if we 
consider connected DGLs over $\Q$, we obtain the usual Lie-Massey products
of rational homotopy theory 
(cf.\ \cite{AlldR,RetaL} or \cite[V.1]{TanrH}).

\subsection{Generalized Massey products}
\label{sgmp}\stepcounter{thm}

If \ $K=(\underline{K}_{n})_{n=0}^{\infty}$ \ is any associative ring 
spectrum, one can define the corresponding \emph{Massey products} by
dualizing diagram \ref{fig14} above, replacing the Whitehead products \ 
$w_{r,s}:S^{r+s-1}\to S^{r}\vee S^{s}$ \ by the multiplication maps \ 
$m^{r,s}:\underline{K}_{r}\times \underline{K}_{s}\to \underline{K}_{r+s}$:

However, in this case it is usual to dualize our description of the
higher order homotopy operation associated to the corresponding 
homotopy-commutative diagram \ $\A:\Gamma\to\Ta$, \ and defining \ 
$\llrra{\A}$ \ to be the collection of adjoints in \ 
$[X,\Omega\underline{K}_{r+s+t}]$ \ to the classes \ 
$\llrr{\A}\subseteq [\Sigma X,\underline{K}_{r+s+t}]$ \ 
(compare \cite{GPorH}).

When \ $K=HR$ \ is the Eilenberg-Mac Lane spectrum for the ring $R$, we
obtain the usual Massey products in \ $H^{\ast}(X;R)$ \ (see \cite{MassN} and
\cite[\S 2]{MUehJ}), and their higher-order generalizations 
(cf.\ \cite{KraiM,MilgI}). The $K$-theory and cobordism versions have been 
considered in \cite{SnaiM} and \cite{JCAlexC} respectively, (and elsewhere).
Of course, one can define such products in any model category with
products \ -- \ e.g., for differential graded algebras.


\begin{thebibliography}{ABCDEF}
%
\bibitem[As]{AdHI}
%
J.F.~Adams,
``On the non-existence of elements of {Hopf} invariant one'',\hsm
\textit{Ann.\ Math.\ (2)} \textbf{72} (1960), No.\ 1, pp.\ 20-104.
%
\bibitem[Am]{AdemI}
J.~Adem,
``The iteration of the {Steenrod} squares in algebraic topology'',\hsm
\textit{Proc.\ Nat.\ Acad.\ Sci.\ USA}\textbf{38} (1952),pp.\ 720-726.
%
\bibitem[Ald]{AlldR}
C.~Allday,
``Rational Whitehead products and a spectral sequence of Quillen'',\hsm
\textit{Pac.\ J.\ Math.} \textbf{46} (1973) No.\ 2, pp.\ 313-323.
%
\bibitem[Ale]{JCAlexC}
J.C.~Alexander,
``Cobordism {Massey} products'',\hsm
\textit{Trans. AMS} \textbf{166} (1972), 197-214.
%
\bibitem[BJM]{BJMahT}
M.G.~Barratt, J.D.S.~Jones \& M.E.~Mahowald,
``Relations amongst Toda brackets and the Kervaire invariant in dimension 
$64$'',\hsm
\textit{J.\ Lond.\ Math.\ Soc.} \textbf{30} (1984), pp.\ 533-550.
%
\bibitem[Ba]{BauG}
H.J.~Baues,
``Geometry of loop spaces and the cobar construction'',\hsm
\textit{Mem.\ AMS} \textbf{230}, AMS, Providence, RI, 1980.
%
\bibitem[B1]{BlaH}
D.\ Blanc,
``A Hurewicz spectral sequence for homology'',\hsm
\textit{Trans.\ AMS} \textbf{318} (1990) No.\ 1, pp.\ 335-354.
%
\bibitem[B2]{BlaO}
D.\ Blanc,
``Operations on resolutions and the reverse {Adams} spectral sequence'',\hsm
\textit{Trans.\ AMS} \textbf{342} (1994) No.\ 1, pp.\ 197-213.
%
\bibitem[B3]{BlaHH}
D.\ Blanc,
``Higher homotopy operations and the realizability of homotopy groups'',\hsm 
\textit{Proc.\ Lond.\ Math.\ Soc.\ (3)} \textbf{70} (1995), pp.\ 214-240.
%
\bibitem[B4]{BlaHO}
D.~Blanc,
``Homotopy operations and the obstructions to being an $H$-space'',\hsm
\textit{Manus.\ Math.} \textbf{88} (1995) No.\ 4, pp.\ 497-515.
%
\bibitem[Bo]{BoaH}
J.M.~Boardman,
``Homotopy structures and the language of trees'',\hsm
In \textit{Algebraic Topology}, Proc. Symp. Pure Math. \textbf{22},
AMS, Providence, RI, 1971, pp.\ 37-58.
%
\bibitem[BV]{BVogHI}
J.M.~Boardman \& R.M.~Vogt,
\textit{Homotopy Invariant Algebraic Structures on Topological Spaces},\hsm
Springer-\-Verlag \textit{Lec.\ Notes Math.} \textbf{347}, 
Berlin-\-New York, 1973.
%
\bibitem[CP]{CPorV}
J.-M.~Cordier \& T.~Porter,
``Vogt's theorem on categories of homotopy coherent diagrams'',\hsm
\textit{Math. Proc. Camb. Phil. Soc.} \textbf{100} (1986), No. 1, pp.\ 65-90.
%
\bibitem[DKS]{DKSmH}
W.G.~Dwyer, D.M.~Kan, \& J.H.~Smith,
``Homotopy commutative diagrams and their realizations'',\hsm
\textit{J. Pure Appl. Alg.}, \textbf{57} (1989), No. 1, pp.\ 5-24.
%
\bibitem[H]{HoltH}
D.N.~Holtzman,
``Higher order cohomology operations in the {$p$}-torsion-free category'',\hsm
\textit{Neder. Akad. Weten. Proc.} \textbf{44} (1982), No.\ 2, pp.\ 183-200.
%
\bibitem[Kl]{KlauT}
S.~Klaus, 
``Towers and Pyramids,  I'',
\textit{Fund.\ Math} \textbf{13} (2001), No. 5, pp.\ 663-683.
%
\bibitem[KK]{KKrisS}
A.~Kock \& L.~Kristensen,
``A secondary product stucture in cohomology theory'',\hsm
\textit{Math. Scand.} \textbf{17} (1965), pp.\ 113-149.
%
\bibitem[K]{KraiM}
D.P.~Kraines,
``Massey higher products'',\hsm
\textit{Trans. AMS} \textbf{124} (1966), 431-449.
%
\bibitem[Kr]{KristS}
L.~Kristensen,
``On secondary cohomology operations'',
\textit{Math. Scand.} \textbf{12} (1963), pp.\ 57-82.
%
\bibitem[Mc]{MacC}
S.\ {Mac Lane},
\textit{Categories for the Working Mathematician},\hsm
Springer-\-Verlag \textit{Grad.\ Texts in Math.} \textbf{5}, 
Berlin-\-New York, 1971.
%
\bibitem[MP]{MPetS}
M.E.~Mahowald \& F.P.~Peterson,
``Secondary operations on the {Thom} class'',\hsm
\textit{Topology} \textbf{2} (1964), pp.\ 367-377
%
\bibitem[Ms]{MassN}
W.S.~Massey,
``A new cohomology invariant of topological spaces'',\hsm
\textit{Bull.\ AMS} \textbf{57} (1951), p.\ 74.
%
\bibitem[MU]{MUehJ}
W.S.~Massey \& H.~Uehara,
``The Jacobi identity for Whitehead products'',\hsm
in \textit{Algebraic geometry and topology}, Princeton U. Press, Princeton, 
1957, pp.\ 361-377.
%
\bibitem[Mar]{MargS}
H.R.~Margolis,
\textit{Spectra and the Steenrod Algebra: \ Modules over the Steenrod Algebra 
and the Stable Homotopy Category},\hsm
North-Holland, Amsterdam-\-New York, 1983.
%
\bibitem[Mau]{MaunC}
C.R.F.~Maunder,
``Cohomology operations of the {$N$}-th kind'',\hsm
\textit{Proc. Lond. Math. Soc. Ser. (2)} \textbf{13} (1963), pp.\ 125-154.
%
\bibitem[May]{MayG}
J.P.~May,
\textit{The Geometry of Iterated Loop Spaces},\hsm
Springer-\-Verlag \textit{Lec.\ Notes Math.} \textbf{271}, 
Berlin-\-New York, 1972.
%
\bibitem[Mi]{MilgI}
R.J.~Milgram,
``Iterated loop spaces'',\hsm
\textit{Ann.\ Math.\ (2)} \textbf{84} (1966), pp.\ 386-403.
%
\bibitem[Mo]{MoriHT}
M.~Mori,
``On higher {Toda} brackets'',\hsm
\textit{Bull. College Sci. Univ. Ryukyus} \textbf{35} (1983), pp.\ 1-4.
%
\bibitem[PS]{PSteS}
F.P.~Peterson \& N.~Stein,
``Secondary cohomology operations: two formulas'',\hsm
\textit{Amer.\ J.\ Math.} \textbf{81} (1959), pp.\ 231-305.
%
\bibitem[P1]{GPorW}
G.J.~Porter,
``Higher order {Whitehead} products'',\hsm
\textit{Topology} \textbf{3} (1965), 123-165.
%
\bibitem[P2]{GPorH}
G.J.~Porter,
``Higher products'',\hsm
\textit{Trans.\ AMS} \textbf{148} (1970), 315-345.
%
\bibitem[R]{RetaL}
V.S.~Retakh,
``Lie-Massey brackets and $n$-homotopically multiplicative maps of
differential graded Lie algebras'',\hsm
\textit{J.\ Pure \& Appl.\ Alg.} \textbf{89} (1993) No.\ 1-2, pp.\ 217-229.
%
\bibitem[Sc]{SchouA}
P.-H.\ Schoute,
``Analytic treatment of the polytopes regularly derived from the regular
polytopes'',\hsm 
\textit{Verh.\ Kon.\ Akad.\ Wet.\ Amst.}, \textbf{11}, 1911.
%
\bibitem[SV]{SchVo}
R.\ Schw\"anzl \& R.M.\ Vogt,
``Coherence in homotopy groups actions'',
in {\it Transformation Groups, Poznan 1985}, 
Lecture Notes in Mathematics 1217 (1986), pp.\ 363--390.                
%
\bibitem[Se]{SegCC}
G.B.~Segal,
``Categories and cohomology theories'',\hsm
\textit{Topology} \textbf{13} (1974), pp.\ 293-312.
%
\bibitem[SS]{SSterQ}
S.~Shnider \& S. Sternberg,
\textit{Quantum groups: from coalgebras to Drinfel'd algebras},\hsm
International Press \textit{Grad.\ Texts in Math.\ Phys.} \textbf{II}, 
Cambridge, MA, 1993.
%
\bibitem[Sn]{SnaiM}
V.P.~Snaith,
``Massey products in {$K$}-theory'',\hsm
\textit{Proc. Camb. Phil. Soc.} \textbf{68} (1970), 303-320.
%
\bibitem[Sp1]{SpanS}
E.H.~Spanier,
``Secondary operations on mappings and cohomology'',\hsm
\textit{Ann.\ Math.\ (2)} \textbf{75} (1962) No.\ 2, pp.\ 260-282.
%
\bibitem[Sp2]{SpanH}
E.H.~Spanier,
``Higher order operations'',\hsm
\textit{Trans.\ AMS} \textbf{109} (1963), pp.\ 509-539.
%
\bibitem[St1]{StaH1}
J.D.\ Stasheff,
``Homotopy associativity of {$H$}-spaces, {I}'',\hsm
\textit{Trans.\ AMS} \textbf{108} (1963) pp.\ 275-292.
%
\bibitem[St2]{StaH2}
J.D.\ Stasheff,
``Homotopy associativity of {$H$}-spaces, {II}'',\hsm
\textit{Trans.\ AMS} \textbf{108} (1963) pp.\ 293-312.
%
\bibitem[St3]{StaH}
J.D.\ Stasheff,
\textit{$H$-spaces from a Homotopy Point of View},\hsm
Springer-\-Verlag \textit{Lec.\ Notes Math.} \textbf{161}, 
Berlin-\-New York, 1970.
%
\bibitem[Ta]{TanrH}
D.~Tanr\'{e},
\textit{Homotopie Rationelle: Mod\`{e}les de Chen, Quillen, Sullivan},\hsm
Springer-\-Verlag \textit{Lec.\ Notes Math.} \textbf{1025}, Berlin-\-New
York, 1983.
%
\bibitem[T1]{TodG}
H.\  Toda,
``Generalized {Whitehead} products and homotopy groups of spheres'',\hsm
\textit{J.\ Inst.\ Polytech.\ Osaka City U., Ser.\ A, Math.} \textbf{3} 
(1952), pp.\ 43-82.
%
\bibitem[T2]{TodC}
H.\ Toda,
\textit{Composition methods in the homotopy groups of spheres},\hsm
Adv.\ in Math.\ Study \textbf{49}, Princeton U. Press, Princeton, 1962.
%
\bibitem[To]{TonkR}
A.P.~Tonks,
``Relating the asso\-cia\-hed\-ron and the permu\-to\-hed\-ron'',\hsm
in J.L.\ Loday, J.D.\ Stasheff, \& A.A.\ Voronov, eds., 
\textit{Operads: Proceedings of Renaissance Conferences 
(Hartford,CT/Luminy, 1995} Contemp.\ Math. \textbf{202}, AMS, Providence, 
RI 1997, pp.\ 33-36.
%
\bibitem[V]{VogtC}
R.M.~Vogt,
``Cofibrant operads and universal {$E_\infty$} operads'',\hsm
in R.M.~Vogt editor, \textit{Workshop on Operads, Osnabr\"uck, June
1998}, Preprint of Univ. Bielefeld 1999, pp.~81-89.

\bibitem[Wa]{GWalkL}
G.\ Walker,
``Long Toda brackets'',\hsm
in \textit{Proc.\ Adv.\ Studies Inst.\ on Algebraic Topology, vol.\ III},
Aarhus U.\ Mat.\ Inst.\ Various Publ.\ Ser.\ \textbf{13}, Aarhus 1970,
pp.\ 612-631.
%
\bibitem[Wh]{GWhHT}
G.W.~Whitehead,
\textit{Homotopy Theory},\hsm
M.I.T. Press, Cambridge, MA, 1953.
%
\end{thebibliography}
\end{document}